\documentclass[12pt]{amsart}
\usepackage{amsmath,amssymb}
\usepackage{graphicx}
\title{Geometry of lines on certain Moishezon threefolds. I.
Explicit description of families of twistor lines}
\author{Nobuhiro Honda$^{\dag}$
}
\thanks
{$^{\dag}$This work was partially supported by
Research Fellowships of the 
Japan Society for the Promotion
of Science for Young Scientists.
}

\newcommand{\ol}{\overline}
\newcommand{\ra}{\rightarrow}
\newcommand{\lra}{\longrightarrow}
\newcommand{\da}{\downarrow}
\newcommand{\ua}{\uparrow}

\newcommand{\set}{\,|\,}
\newcommand{\proofend}{\hfill$\square$}

\newtheorem{prop}{Proposition}[section]
\newtheorem{lemma}[prop]{Lemma}

\newtheorem{thm}[prop]{Theorem}

\newtheorem{cor}[prop]{Corollary}

\newtheorem{definition}[prop]{Definition}

\begin{document}

\maketitle

\begin{abstract} We study real lines  on certain Moishezon threefolds
which are potentially   twistor spaces of $3\mathbf{CP}^2$.  Here,  
line means a smooth rational curve whose normal bundle is $
O(1)^{\oplus 2}$ and the reality implies the invariance under an
anti-holomorphic involution on  the threefolds. 
Our threefolds are birational to double covering of $\mathbf{CP}^3$
branched along a singular quartic surface.
On these threefolds we  find  families of real lines
in explicit form and prove  that which families have to be
 chosen as  twistor lines depend on how we take small resolutions
 of the double covering. 
 This is a first step
for determining the moduli space of self-dual metrics on
$3\mathbf{CP}^2$ of positive scalar curvature,  which admit  a
non-trivial Killing field. 
\end{abstract}

\tableofcontents


\section{Introduction} A Riemannian metric on an oriented four-manifold
is called self-dual iff the anti-self-dual part of the Weyl conformal
curvature of the metric identically  vanishes. Basic examples are
provided by the round metric on the four-sphere and the Fubini-Study
metric on the  complex projective plane. In general, one can expect that
if two four-manifolds admit self-dual metrics respectively, then their
connected sum  will also admit a self-dual metric. In fact, Y.S. Poon
\cite{P86}  constructed  explicit examples of  self-dual metrics on
$2\mathbf{CP}^2$, the connected sum of two complex projective planes.
He further showed that on $2\mathbf{CP}^2$  there are no self-dual
metrics other than his metrics, under the assumption of 
the positivity of the scalar curvature.  Later, C. LeBrun \cite{LB91}
and D. Joyce \cite{J95} respectively constructed self-dual metrics of
positive scalar curvature on $n\mathbf{CP}^2$ for any $n\geq 1$. These
are  called LeBrun metrics and Joyce metrics respectively, and
have nice characterizations by the (conformal) isometry group. Namely, 
A. Fujiki \cite{F00} proved that if a self-dual metric on
$n\mathbf{CP}^2$  has $U(1)\times U(1)$ as the identity component of the
isometry group, then it must be a Joyce metric. LeBrun \cite{LB93} 
showed that  if a self-dual metric on $n\mathbf{CP}^2$ has a non-trivial
semi-free
$U(1)$-isometry, then the metric must be a LeBrun metric. Here, a
$U(1)$-action on a manifold $M$ is called semi-free if the isotropy
group is
$U(1)$ or identity only, at every point of $M$.

Now drop the assumption of  the semi-freeness for $U(1)$-isometry. 
We further suppose that the metric is not a Joyce metric,
since it is immediate to get non-semi-free $U(1)$-isometries by choosing 
subgroups of $U(1)\times U(1)$.
Then in \cite{Hon02,Hon-II}
the author has   shown  that, for $n=3$, there already 
exist  self-dual metrics of positive scalar curvature
whose identity components of the
isometry groups are $U(1)$ acting non-semi-freely on $3\mathbf{CP}^2$.
 (Note that Poon's metrics on 
$2\mathbf{CP}^2$ coincide with Joyce metrics and hence there is no self-dual metric
on $2\mathbf{CP}^2$
whose isometry group is $U(1)$.)
 The purpose of the present paper is to study the twistor spaces of these
self-dual metrics on $3\mathbf{CP}^2$. In other words, we are interested in twistor
spaces of
$3\mathbf{CP}^2$  which has a non-trivial holomorphic $U(1)$-action but which are
different from LeBrun  twistor spaces.
Note that on $3\mathbf{CP}^2$ Joyce metrics coincide with LeBrun
metrics with torus action.

By a result of Kreu\ss ler-Kurke \cite{KK92}  such
a twistor space is always Moishezon.
First in Section \ref{s-defeq} we determine the defining equations of 
projective models of the above mentioned  twistor spaces of
$3\mathbf{CP}^2$ (Propositions \ref{prop-def-B} and \ref{prop-necessa}). 
Here note that by a famous theorem of Hitchin \cite{Hi81} the twistor
spaces themselves  cannot be projective algebraic, so the projective
models are only birational to  the real  twistor spaces.
Roughly the projective models have a structure of  the double covering
of $\mathbf{CP}^3$ branched along a singular quartic surface which is 
birational to an elliptic ruled surface.
The defining equations allow us to prove  that {\em our non-semi-free
$U(1)$-action on $3\mathbf{CP}^2$ is uniquely determined up to
equivariant diffeomorphims} (Proposition \ref{prop-class}).

Then it naturally arises a  question  whether, conversely, the
singular projective threefold having this kind of structure  always becomes
birational to  a twistor space. Also we want to know how one can obtain the twistor
spaces from the projective threefolds by means of resolution of singularities.
These are fascinating but  difficult problems, since they are equivalent to
determine the moduli space of self-dual metrics on
$3\mathbf{CP}^2$ of positive scalar curvature, which admit a  non-trivial
Killing field.

The most  primitive way to
 prove that given threefold is actually a twistor space is  to find a 
family of  twistor lines. Here a twistor line is by definition a fiber
of the twistor fibration
$Z\ra M$, where $M$ is the base four-manifold. In Section \ref{s-detc},
we see that the image of a real line 
in our threefold  under the double covering map onto
$\mathbf{CP}^3$ is a line iff the real line  intersects  some real smooth
rational curve, 
and that otherwise the image is a real
conic  whose  intersection number with the branch quartic surface is
at least two for any  intersection points
(Proposition \ref{prop-image}). Following I. Hadan \cite{Ha00}, we will call the
latter kind of curve {\em{a touching conic}}. 

Thus in order to find all the twistor lines, we need to study the space
of  real touching conics for our branch quartic surface. 
Note that every conic in
$\mathbf{CP}^3$ is contained in a plane. In the rest of Section
\ref{s-detc} we explicitly write down the defining equations of real touching
conics contained in real
$\mathbf C^*$-invariant planes, where such planes are parametrized by a circle.  
Our results show that for each of such a real plane, real touching conics in it
form just two  real one-dimensional families (Propositions
\ref{prop-a},
\ref{prop-b} and
\ref{prop-c} for precisely.) These families are not  pencil in general.

In Section \ref{s-inv} we  study the inverse images of real touching conics
classified in the previous section. The inverse images are of course candidates of
twistor lines. It is immediate to see that the inverse image
has just two irreducible components which are mapped biholomorphically onto
the conic.  
We show that these two components form the real parts of   real pencils  on a real
$U(1)$-invariant smooth surface (which is the inverse image of a real $U(1)$-invariant
plane).

Next in Section \ref{s-nb} we  calculate the normal bundle of these two
irreducible components inside a smooth threefold which is a small resolution of the
singular projective model (the double cover of
$\mathbf{CP}^3$). 
It is easy to show that the normal bundle is always isomorphic to
 $O(1)^{\oplus
2}$ or $O\oplus O(2)$. But it is quite subtle  which 
one actually occurs.
 In fact, we will see that {\em the result depends on the choice of  small
resolutions}, and  that {\em which one of the two irreducible components 
must be chosen
for twistor line depends on which resolution we take}:
for each  small resolution, we define a function defined on the parameter space of
 real invariant planes,
 and show that the normal bundle degenerates into $O\oplus O(2)$ precisely when
the function has a critical point. We then  determine the critical points
of the functions. As a result, we 
 can exactly determine the images of twistor lines contained in the real $\mathbf
C^*$-invariant planes (Theorem \ref{thm-type}).
Furthermore, these calculations enable us to prove  that {\em among 
twenty-four ways of possible small resolutions, only two resolutions
can yield a twistor space} (Theorem \ref{thm-sr}), and also to  determine which
irreducible component must be chosen, for each of the two resolutions.

As a byproduct of these investigations,  we can
prove that {\em in our smooth threefolds there always exist  families
of  real smooth
rational curves whose normal bundles are
$O(1)^{\oplus 2}$ but  which cannot be  twistor lines, by showing that 
they can be deformed into a smooth rational 
curve whose normal bundle is
$O\oplus O(2)$ while keeping the reality} (Corollary \ref{cor-chatl} and
Proposition \ref{prop-nottl}). Thus  
the twistor line is not characterized by this
property in general. However, it should be noted that  since the Penrose
correspondence between  twistor spaces and self-dual metrics is  purely local matter,
these `fake twistor lines' certainly define a self-dual metric on some non-compact
four-manifold.
It will be interesting to study   what the four-manifold is, and  to clarify the
behavior of the self-dual metric near the ends.

Finally in Section \ref{s-line} we calculate the normal bundle of an
irreducible component of the inverse image of a real line in $\mathbf{CP}^3$,
 which is a
candidate of twistor lines mentioned above, and show that it is always
real and has $O(1)^{\oplus 2}$ as the normal bundle.
(Proposition \ref{prop-imline}). Thus the
situation is quite simple in this case. A connection between the subjects in
the previous sections is that these lines in $\mathbf{CP}^3$ can be considered
as a  limit  of touching conics  studied in Section
\ref{s-detc}--\ref{s-nb}. This fact will play an important role  for forthcoming
investigation of the compactification of the space of real lines in our
threefolds.

Finally I mention the relationship between Hadan's 
elaborate work \cite{Ha00} and ours.
Broadly speaking, the threefolds Hadan studied are candidates of the twistor 
spaces of {\em generic} self-dual metrics on $3\mathbf{CP}^2$
so that they have {\em no symmetry},
whereas ours should be the twistor spaces of self-dual metrics on $3\mathbf{CP}^2$
{\em with  semi-free} $U(1)$-{\em symmetry}.
Although both threefolds are birational to the double covering of $\mathbf{CP}^3$
branched along a quartic surface with isolated singularities,
and although (as far as I could understand) he also aimed to determine
all of the twistor lines,
there are quite a few differences.
First of all, the branch quartic surface is birational to a K3 surface
in Hadan's case, while it is birational to an elliptic ruled surface
in our case.
Secondly, since his branch quartic surface has much more singularities
than ours, there are much more singular plane sections of the quartic.
This yields monodromy problem of families of touching conics,
which is the main subject of his investigation.
Thirdly, he mainly concerns lines which are not necessarily
real. But we only consider real ones, and the reality plays a crucial role.
Fourthly, he focused his attention to the behavior of generic  lines,
and therefore did not aim to get the defining equations of touching conics.
In contrast, we are concerning twistor lines which can be considered as
a limit of generic twistor lines.
Finally, his study is mostly on $\mathbf{CP}^3$, and does not 
investigate the inverse images of touching conics.
In particular, he does not consider their normal bundles,
whereas a considerable part of our investigation is 
devoted to determining them.
As a consequences of these differences, there are  few overlaps
between his results and ours.

I would like to thank Nobuyoshi Takahashi for stimulating discussions
and helpful  comments.

\section{Defining equations of the branch quartic surfaces and their
singularities}
\label{s-defeq} Let $g$ be a self-dual metric on $3\mathbf{CP}^2$ of positive
scalar curvature, and assume that $g$ is not conformally isometric to LeBrun
metrics. Let
$Z$  be the twistor space of $g$, and denote by $(-1/2)K_Z$  the fundamental
line bundle which is the canonical square root of the anticanonical line bundle
of $Z$.
These non-LeBrun twistor spaces of $3\mathbf{CP}^2$ are extensively studied in 
 Kreu\ss ler and Kurke  \cite{KK92} and Poon \cite{P92}, 
and it has been proved that  the fundamental
system (the complete linear system associated to the fundamental line bundle) is
free and of three-dimensional, and induces a  surjective morphism
$\Phi:Z\ra\mathbf{CP}^3$ which is generically two-to-one, and  that the branch 
divisor $B$ is a quartic surface with only isolated singularities.
Furthermore, there is the following diagram:

$$
\begin{array}{rlll}
&Z&&\\
\mu&\da&\searrow\Phi&\\
&Z_0&\lra\mathbf{CP}^3\\
&&\Phi_0&
\end{array}
$$
where $\Phi_0:Z_0\ra \mathbf{CP}^3$
denotes the  double covering branched along
$B$, and $\mu$ is {\it{a small
resolution}}\, of the singularities of $Z_0$ over the singular points of
$B$. 

Generically, $B$ has only ordinary double points \cite{KK92, P92} and
hence is birational to a K3 surface. As a consequence, 
one can deduce that $Z$ does not
admit a non-zero holomorphic vector field.
However, the author showed
in \cite{Hon02} and \cite{Hon-II}  that 
if $B$ degenerates to have non-ADE singularities, then 
$Z$ admits a non-zero holomorphic vector field,
and that such a twistor space of $3\mathbf{CP}^2$ actually exists.
Concerning a defining equation of
the branch quartic $B$ for such twistor spaces, we have the following
proposition which is the starting point of our investigation.

\begin{prop}
\label{prop-def-B} Let $g$ be a non-LeBrun self-dual metric on 
$\mathbf{CP}^2$ of positive scalar curvature, and assume the
existence of a non-trivial Killing field. Let 
$\Phi:Z\ra\mathbf{CP}^3$ and $B\subset\mathbf{CP}^3$ be as above.
Then there exists a homogeneous coordinate
$(y_0:y_1:y_2:y_3)$ on
$\mathbf{CP}^3$ fulfilling (i)-(iii) below: 

\noindent (i) a defining equation of
$B$ is given  by
\begin{equation}\label{eqn-B}
\left(y_2y_3+Q(y_0,y_1)\right)^2-y_0y_1(y_0+y_1)(ay_0-by_1)=0,
\end{equation} where $Q(y_0,y_1)$ is a   
quadratic form of $y_0$ and
$y_1$ with real coefficients, and $a$ and $b$ are  positive  real numbers,

\noindent (ii) the naturally induced real structure on $\mathbf{CP}^3$ is given by
$$
\sigma(y_0:y_1:y_2:y_3)=\left(\ol{y}_0:\ol{y}_1:\ol{y}_3:\ol{y}_2\right),
$$

\noindent (iii) the naturally induced $U(1)$-action on
$\mathbf{CP}^3$ is given by
$$(y_0:y_1:y_2:y_3)\mapsto 
\left(y_0:y_1:e^{i\theta}y_2:e^{-i\theta}y_3\right),\hspace{3mm} 
e^{i\theta}\in U(1).$$
\end{prop}
 
\noindent Proof.  If the fundamental system of $Z$ 
is free, there are just four reducible members, all of which are real
\cite{P92,KK92}. 
We write $\Phi^{-1}(H_{ i})=D_{ i}+\ol{D}_{ i}$, $1\leq i\leq 4$, where
$H_i$ is a real plane in $\mathbf{CP}^3$.
 The restrictions of $\Phi$ onto $D_{i}$ and
$\ol{D}_{i}$ are obviously birational morphisms onto $H_{i}$, so that,
together with the fact that $(-1/2)K_Z\cdot L_i=2$,  it can be readily seen that
$C_{i}:=\Phi(L_{i})$ is a conic contained in $B$. This implies that
the restriction of 
$B$ onto $H_{i}$ is a  conic of multiplicity two. Namely, $C_{i}$,
$1\leq i\leq 4$, is so called {\em a trope} of $B$. 

A Killing field naturally gives rise to an isometric $U(1)$-action,
which can be canonically lifted to a holomorphic $U(1)$-action on the
twistor space.
This action  naturally goes down to 
$\mathbf{CP}^3$, and every subvarieties above are clearly preserved by
these $U(1)$-actions. In particular, $C_{i}$ is a
$U(1)$-invariant conic on a
$U(1)$-invariant plane $H_{i}$, where the $U(1)$-action
is induced by the vector field. Since the twistor fibration is 
$U(1)$-equivariant and generically one-to-one on $D_i$, and since $\Phi|_{D_i}:D_i\ra
H_i$ is also $U(1)$-equivariant and birational, the  $U(1)$-action on any
$H_{i}$ is non-trivial. Hence  $U(1)$ acts non-trivially on $C_{i}$. For
$j\neq i$, put
$l_{ij}=H_i\cap H_j$, which is clearly a real $U(1)$-invariant  line. Then 
 $C_{i}\cap
l_{ij}$ must be the two
$U(1)$-fixed points on $C_{i}$, since it is real set and since there is 
no real point on $C_i$. This implies that $l_{ij}$ is
independent of the choice of $j\neq i$. So we write
$l_{ij}=l_{\infty}$ and let $P_{\infty}$ and $\ol{P}_{\infty}$ be  the
two fixed points of the $U(1)$-action on $l_{\infty}$. Then 
$H_{i}$,
$1\leq i\leq 4$, must be  real  members of the real pencil of planes whose base locus is
$l_{\infty}$.  Since
$l_{\infty}$  is a real line, we can choose  real linear forms $y_0$ and
$y_1$ such that
$l_{\infty}=\{y_0=y_1=0\}$. Further, since any of $H_{i}$ is real, by
applying a real projective transformation (with respect to
$(y_0,y_1)$), we may assume that $\cup_{i=1}^4
H_{i}=\{y_0y_1(y_0+y_1)(ay_0-by_1)=0\}$, where $a,b\in\mathbf R$ with
$a>0$ and
$b>0$. 

 As seen above, every $C_{i}$ goes through $P_{\infty}$ and
$\ol{P}_{\infty}$.
 Let $l_{i}$, $1\leq i\leq 4$  be the tangent line of
$C_{i}$ at $P_{\infty}$. Now we claim that $l_{i}$'s  are lying on the
same plane. Let
$H$ be the plane containing $l_1$ and $l_2$.  
Then by using 
$l_1\cap B=P_{\infty}= l_2\cap B$,
we can easily deduce that  $B\cap H$  is a 
 union of lines, all of which goes through $P_{\infty}$. 
Suppose that $l_3$ is not contained in $H$. Then the line $H\cap
H_3$ is not tangent to $C_3$, so there is an intersection point  of
$C_3\cap H$ other than
$P_{\infty}$. Then the line $H\cap H_3$ is contained in $B$, because we
have already seen that $B\cap H$ is a union of lines all passing through
$P_{\infty}$. This is a contradiction since
$B\cap H_3=2C_3$. Similarly we have
$l_4\subset H$. Therefore $l_{i}\subset H$ for any $i$, as claimed.  Because
$C_{i}$'s are real, the plane
$\sigma(H)$  contains the tangent lines of
$C_{i}$'s at $\ol{P}_{\infty}$. Let $y_3$ be a linear form on
$\mathbf{CP}^3$ defining $H$, and set
$y_2:=\ol{\sigma^*y_3}$. Then  $(y_0:y_1:y_2:y_3)$ is a homogeneous
coordinate on $\mathbf{CP}^3$. By our choice, we have
$\sigma(y_0:y_1:y_2:y_3)=(\ol{y}_0:\ol{y}_1:\ol{y}_3:\ol{y}_2)$,
$P_{\infty}=(0:0:0:1)$ and $\ol{P}_{\infty}=(0:0:1:0)$.

Because the planes $\{y_i=0\}$ are $U(1)$-invariant, our $U(1)$-action can be linearized with respect to the homogeneous coordinate
$(y_0:y_1:y_2:y_3)$. 
Further, since $H_i$'s are $U(1)$-invariant, the action can be written
$(y_0:y_1:y_2:y_3)\mapsto (y_0:y_1:e^{i\alpha\theta}y_2:
e^{i\beta\theta}y_3)$ for $e^{i\theta}\in U(1)$, where $\alpha$ and $\beta$ are
relatively prime integers.
Moreover, since the conics $C_i$'s are $U(1)$-invariant, we can suppose
$\alpha=1$, $\beta=-1$. Thus the $U(1)$-action can be written in the form
(ii) of the proposition.

Let $F=F(y_0,y_1,y_2,y_3)$ be a defining equation of $B$. Since $B$ is
$U(1)$-invariant,  monomials appeared  in $F$ must be  in the ideal
$(y_2y_3, y_0^2,y_0y_1,y_1^2)^2$.
Moreover, 
$F$ contains the monomial $y_2^2y_3^2$, since otherwise the restriction onto
$\{y_1=0\}$ would be the union of two different conics, 
which contradict to the fact that $C_i$ is a trope.
We assume that its coefficient  is $1$. Then
$F$ can be written in the form $(y_2y_3+Q(y_0,y_1))^2-q(y_0,y_1)$, 
where
$Q(y_0,y_1)\in (y_0,y_1)^2$  and 
$q(y_0,y_1)\in (y_0,y_1)^4$  are uniquely determined polynomials
with real coefficients. Then it again follows from $C_i$ being a trope that 
$q(y_0,y_1)=ky_0y_1(y_0+y_1)(ay_0-by_1)$ for some constant $k\in\mathbf
R^{\times}$. 
Finally, if $k$  is negative, exchange $y_0$ and $y_1$. 
Then we get the equation of the form (\ref{eqn-B}), 
and we have proved all of the claims in the proposition.
\proofend

\vspace{2mm}


Next we study the singular locus of $B$.

\begin{prop}\label{sing. of B} Let $B$ be a real quartic surface defined by
the equation
$$\left(y_2y_3+Q(y_0,y_1)\right)^2-y_0y_1(y_0+y_1)(ay_0-by_1)=0$$ where
$Q(y_0,y_1)$ is a real quadratic form of $y_0$ and $y_1$,
and
$a$ and $b$ are  real positive numbers.
Let
$A$ be the set
$\{(y_0:y_1:0:0)\set (y_0:y_1)$ is a multiple root of the quartic
equation
$Q(y_0,y_1)^2-y_0y_1(y_0+y_1)(ay_0-by_1)=0$$\}$.  (Here we think
of this as an equation
on
$\mathbf{CP}^1=\{(y_0:y_1)\}$.) Then we have:
(i) {\em Sing}$(B)=\{P_{\infty},\ol{P}_{\infty}\}\cup A$, where we put
$P_{\infty}=(0:0:0:1)$,
(ii) $P_{\infty}$ and $\ol{P}_{\infty}$ are elliptic singularities of
type
$\tilde{E}_7$, and
(iii) if $Q(y_0,y_1)\neq 0$, then $(y_0:y_1:0:0)\in A$ is an ordinary double point
iff its multiplicity is two.
\end{prop}
 In particular,  every singular point of $B$ is isolated.

\vspace{2mm}
\noindent Proof. (i) First we show that (Sing$B)\cap\{y_3\neq
0\}=\{P_{\infty}\}$, by calculating the Jacobian. Let $x_i=y_i/y_3$ 
($0\leq i\leq 2$) be affine
coordinates on $y_3\neq 0$. Then the equation of $B$ becomes 
$(x_2+Q(x_0,x_1))^2-x_0x_1(x_0+x_1)(ax_0-bx_1)=0$. Differentiating with
respect to 
$x_2$, we get $x_2+Q(x_0,x_1)=0$ so that we have
$x_0x_1(x_0+x_1)(ax_0-bx_1)=0$. Next differentiating with respect to
$x_0$ and
$x_1$ and then substituting
$x_2+Q(x_0,x_1)=0$, we get
$x_1(x_0+x_1)(ax_0-bx_1)+x_0x_1(ax_0-bx_1)+ax_0x_1(x_0+x_1)=0$ and
$x_0(x_0+x_1)(ax_0-bx_1)+x_0x_1(ax_0-bx_1)-bx_0x_1(x_0+x_1)=0$. From the
former, we obtain that  $x_0=0$ implies $x_1=0$. Then by
$x_2+Q(x_0,x_1)=0$ we have $x_2=0$. Similar argument shows that  if
$x_1, x_0+x_1$ or $ax_0-bx_1$ is zero, then $x_0=x_1=x_2=0$. Conversely,
it is immediate to see that
$(x_0,x_1,x_2)=(0,0,0)$ is a double point of $B$. Thus we get
(Sing$B)\cap\{y_3\neq 0\}=\{P_{\infty}\}$. Because the given
homogeneous polynomial is symmetric with respect to $y_2$ and $y_3$, we have
(Sing$B)\cap\{y_2\neq 0\}=\{\ol{P}_{\infty}\}$.

Next we show that (Sing$B)\cap \{y_2=y_3=0\}=A$. We may suppose $y_1\neq 0$.
Putting
$v_i=y_i/y_1$ for $i=0,2,3$, the equation  of $B$ becomes
$(v_2v_3+Q(v_0,1))^2-f(v_0)=0$, where we put
$f(v_0)=v_0(v_0+1)(av_0-b)$. Substituting $v_2=v_3=0$, we get
$Q(v_0,1)^2-f(v_0)=0$. On the other hand, differentiating with respect to $v_0$ and
substituting $v_2=v_3=0$, we get $(Q(v_0,1)^2-f(v_0))'=0$, where the
prime denotes differential with respect to $v_0$. 
Thus, if $(v_0,v_2,v_3)=(\lambda_0,0,0)$ is a singular point of $B$,
$\lambda_0$ is a multiple root of $Q(v_0,1)^2-f(v_0)=0$.
 Conversely, it is easy to see that
$(v_0,v_2,v_3)=(\lambda_0,0,0)$ is a singular point for  such $\lambda_0$. Thus we get
the claim of (i).

(ii)  is obvious if one notes that we can  use
$(x_0,x_1,x_2+Q(x_0,x_1))$ instead of $(x_0,x_1,x_2)$ as a local
coordinate around $P_{\infty}$.

Finally we show (iii) by using the coordinate $(v_0,v_2,v_3)$ above. Let
$\lambda_0$ be a multiple root of $Q(v_0,1)^2-f(v_0)=0$. Then our
equation of $B$ can be written
$(v_2v_3+2Q(v_0,1))v_2v_3+g(v_0)(v_0-\lambda_0)^2=0$, where $g(v_0)$
is a  polynomial of degree two.
Clearly $\lambda_0$ is a double root iff $g(\lambda_0)\neq 0$.
Suppose $g(\lambda_0)\neq 0$ and define
\begin{equation}\label{coord}
\begin{array}{l}
w_1=\sqrt{g(v_0)}\cdot(v_0-\lambda_0),\\
w_2=\sqrt{2Q(v_0,1)+v_2v_3}\cdot v_2,\\
w_3= \sqrt{2Q(v_0,1)+v_2v_3}\cdot v_3.
\end{array}
\end{equation}
Because $g(\lambda_0)\neq 0$ and  $Q(\lambda_0,1)\neq0$,
$(w_1,w_2,w_3)$ is a local coordinate around $(\lambda_0,0,0)$.
 Then our equation of $B$ becomes
$w_2w_3+w_1^2=0$. Thus the singularity is an ordinary double point.
Conversely, if $g(\lambda_0)=0$, it is immediate to see that our equation of $B$
can be written of the form $w_2w_3+w_1^3=0$ or $w_2w_3+w_1^4=0$ depending on
whether the multiplicity of $\lambda_0$ is three or four.
This implies that $(\lambda_0,0,0)$ is not an ordinary double point.
\proofend

\begin{prop}\label{prop-pos} Let $B$ be as in Proposition \ref{sing. of
B}. Put $f(\lambda)=\lambda(\lambda+1)(a\lambda-b)$. Let
$Z_0\ra\mathbf{CP}^3$ be the double covering branched along $B$. Then if
$Z_0$ admits a small resolution  $Z\ra Z_0$ such that $Z$ is a twistor
space of
$3\mathbf{CP}^2$, then 
$Q(\lambda,1)^2-f(\lambda)\geq 0$ for any $\lambda\in\mathbf R$ and the
equality holds for a unique $\lambda_0\in\mathbf R$. Further, in this case,
the multiplicity of $\lambda_0$ is two.
\end{prop}

\noindent Note that it follows from this proposition that  
$f(\lambda_0)>0$ holds, because we have $Q(\lambda_0,1)^2=f(\lambda_0)$ and
$Q(\lambda_0,1)$ is a real number which is non-zero because otherwise
the restriction $B|_{H_{\lambda_0}}$ would be  $y_0-\lambda_0 y_1=(y_2y_3)^2=0$ that
yields another reducible fundamental divisor.

\vspace{2mm}
\noindent Proof. By  results of Kreu\ss ler \cite{Kr89} and  Kreu\ss
ler-Kurke
\cite{KK92},  we have $\sum(\mu(x)+c(x))=26$ for $Z$ to be a twistor
space of
$3\mathbf{CP}^2$ for a topological reason, where
$\mu(x)$ is the Milnor number of the  singularity $x$ of $B$ and $c(x)$
is the number of irreducible components of a small resolution $Z\ra
Z_0$.  Because elliptic singularity of type
$\tilde{E}_7$ has $\mu=9$ and $c=3$, we get $\sum(\mu(x)+c(x))=2$ for 
other remaining singularities. This implies that there is only one
singularity remaining, and that it must be an ordinary double  point, 
which will be denoted by $P_0$.
Therefore, by
Proposition
\ref{sing. of B} (i), we have
$A=\{P_0\}$. Namely, $Q^2-f=0$ has a unique  multiple root $\lambda_0$.
 The multiplicity
is two by Proposition \ref{sing. of B} (iii).  It is obvious from the uniqueness that
this ordinary double point is  real. Namely, $\lambda_0$ is real.

Next we show that other  solutions of $Q(\lambda,1)^2-f(\lambda)=0$
are not real.  Assume $\lambda\in\mathbf R$, $\lambda\neq\lambda_0$ is a
solution. Then by restricting $B$ to the plane $y_0=\lambda y_1$, we get
$(y_2y_3+Q(\lambda,1)y_1^2)^2-f(\lambda)y_1^4=
(y_2y_3)^2+2Q(\lambda,1)y_2y_3y_1^2+(Q(\lambda,1)^2-f(\lambda))y_1^4=
y_2y_3(y_2y_3+2Q(\lambda,1)y_1^2)$ $(=0)$. Therefore, the point
$(\lambda:1:0:0)$ is a real point of $B$. Since the multiplicity of the
solution $\lambda$ is one, 
 Proposition \ref{sing. of B} shows that this is a smooth point of
$B$. This implies that $Z$ has a real point, contradicting to the
absence of real points on any twistor spaces. Hence the equation
$Q(\lambda,1)^2-f(\lambda)=0$ has no real solution other than
$\lambda_0$. Because $\lambda_0$ is a solution whose multiplicity is
 two, this implies that the polynomial
$Q(\lambda,1)^2-f(\lambda)$ has constant sign on
$\mathbf R\backslash \{\lambda_0\}$.  This sign must be clearly positive.
\proofend

\vspace{2mm}
To investigate the real locus of $B$,  we need the
following elementary 
\begin{lemma}\label{lemma-element} Let $C_{\alpha}=\{y_2y_3=\alpha
y_1^2\}$,
$\alpha\in\mathbf R$ be a real conic in
$\mathbf{CP}^2$, where the real structure is given by
$(y_1:y_2:y_3)\mapsto (\ol{y}_1:\ol{y}_3:\ol{y}_2).$ Then $C_{\alpha}$
has no real point iff
$\alpha<0$.
\end{lemma}

\noindent Proof.
It is immediate to see that the real locus of $C_{\alpha}$ is 
$$\left\{(1:v:\ol{v})\in\mathbf{CP}^2\set |v|=\sqrt{\alpha}\right\}.$$
This is empty iff $\alpha<0$.
\endproof

\begin{prop}\label{prop-realpoint} Let $B$ be as in Proposition
\ref{sing. of B} and suppose that the inequality
$Q(\lambda,1)^2-f(\lambda)\geq 0$ holds on $\mathbf R$ with the
equality holding iff
$\lambda=\lambda_0$ as in Proposition
\ref{prop-pos}. Put $P_0:=(\lambda_0:1:0:0)$, which is clearly a real
point of $B$. Then we have: (i) there is no real point on $B$ other than 
$P_0$ iff the following condition is satisfied:
if $f(\lambda)\geq 0$ and $\lambda\neq\lambda_0$, then
$Q(\lambda,1)>\sqrt{f(\lambda)}$
$\cdots (*)$,
  (ii) if $ (*)$ is satisfied, then there is no real point on any
small resolutions of
$Z_0$.

\end{prop}
\noindent

\vspace{2mm}
\noindent Proof. It is immediate to see that 
 any real point of $B$ is contained in some
real plane
$H_{\lambda}:=\{y_0=\lambda y_1\}$, $\lambda\in\mathbf R\cup\{\infty\}$.
An equation of the restriction
$B_{\lambda}:=B\cap H_{\lambda}$ is given by (as in the proof of
Proposition
\ref{prop-pos}) 
$\left(y_2y_3+Q(\lambda,1)y_1^2\right)^2-f(\lambda)y_1^4=0$. This can be rewritten as
$$B_{\lambda}:
\left\{y_2y_3+\left(Q(\lambda,1)-\sqrt{f(\lambda)}\right)y_1^2\right\}
\left\{y_2y_3+\left(Q(\lambda,1)+\sqrt{f(\lambda)}\right)y_1^2\right\}=0.
$$ Namely, $B_{\lambda}$ is a union of two conics.
If $\lambda\neq\lambda_0$, we have $Q(\lambda,1)^2-f(\lambda)> 0$ by 
our assumption, so both of 
$Q(\lambda,1)-\sqrt{f(\lambda)}$
and $Q(\lambda,1)+\sqrt{f(\lambda)}$ are non-zero.

  Recall that our real structure is given by
$\sigma(y_1:y_2:y_3)=(\ol{y}_1:\ol{y}_3:\ol{y}_2)$ on $H_{\lambda}$
(Proposition \ref{prop-def-B}, (ii)).
Thus each component of $B_{\lambda}$ is real iff the coefficients are
real; namely
$f(\lambda)\geq 0$. Further, the intersection of these two conics is
$\{P_{\infty},\ol{P}_{\infty}\}$. Therefore, there is no real point on
$B_{\lambda}$ if $f(\lambda)<0$. So suppose $f(\lambda)\geq 0$. In this
case, each of the two conics are real, and by Lemma \ref{lemma-element},
both components have no real point iff $Q(\lambda,1)>\sqrt{f(\lambda)}$.
On the other hand, we have $B_{\infty}=\{(y_2y_3+Q(y_0,0))^2=0\}$.
Hence again by Lemma \ref{lemma-element}, we have $Q(1,0)>0$ if
$B_{\infty}$ has no real point. But this   follows from the first
condition.  If $\lambda=\lambda_0$, we have 
$Q(\lambda_0,1)=\sqrt{f(\lambda_0)}$ and hence one of the
components of
$B_{\lambda_0}$ degenerates into a union of two lines whose intersection is $P_0$.
And the other component has no real point since
$Q(\lambda_0,1)+\sqrt{f(\lambda_0)}=2\sqrt{f(\lambda_0)}>0$. Thus we get (i).
 
Next we show that $Z_0$ has no real point other than $\Phi_0^{-1}(P_0)$,
under the condition
$(\ast)$. Suppose $y_1\neq 0$ and use the coordinate $(v_0,v_2,v_3)$
defined in the proof of  Proposition \ref{sing. of B}. Then  $Z_0$ is
given by the equation
$z^2+(v_2v_3+Q(v_0,1))^2-f(v_0)=0$, where $z$ is a fiber coordinate on
the line bundle $ O(2)$ over $\mathbf{CP}^3$.  Thus to prove 
$(Z_0)^{\sigma}=\Phi_0^{-1}(P_0)$ over $y_1\neq 0$,  it suffices to show
that 
\begin{equation}\label{ineq} (v_2v_3+Q(v_0,1))^2-f(v_0)>0
\end{equation} for any real $(v_0, v_2,v_3)\neq (\lambda_0,0,0)$. Recall
that
$(v_0, v_2,v_3)$ is real iff
$v_0\in\mathbf R$ and
$v_2=\ol{v}_3$. Hence (\ref{ineq}) is obvious for real $(v_0, v_2,v_3)$ 
with
$f(v_0)<0$. Assume $f(v_0)\geq 0$. Using the reality condition, we have
$(v_2v_3+Q(v_0,1))^2-f(v_0)=|v_2|^4+2Q(v_0,1)|v_2|^2+
(Q(v_0,1)^2-f(v_0))$.
By our assumption we have $Q(v_0,1)^2-f(v_0)> 0$ for any 
$v_0\in\mathbf R$ with $v_0\neq \lambda_0$. Further, by the condition
$(\ast)$  we have $Q(v_0,1)>\sqrt{f(v_0)}$ for $v_0\neq \lambda_0$
with $f(v_0)\geq 0$. Therefore we have (\ref{ineq}) also for  real
$(v_0, v_2,v_3)\neq (\lambda_0,0,0)$ with $f(v_0)\geq 0$. Thus we have 
$(Z_0)^{\sigma}=\Phi_0^{-1}(P_0)$ over $y_1\neq 0$. 
In the same way we can see that 
$(Z_0)^{\sigma}=\Phi_0^{-1}(P_0)$ over $y_0\neq 0$. 
So it remains to see that there is no real point over the line
$l_{\infty}=\{y_0=y_1=0\}$.
To check this, we introduce a new homogeneous coordinate
$(y_0,y_1,y_2-y_3,y_2+y_3)$ on $\mathbf{CP}^3$.
Then the two subsets $y_2-y_3\neq 0$ and $y_2+y_3\neq 0$ are real,
and it can be easily seen that over the line $l_{\infty}$,
 the equation of $Z_0$ is
of the form $z^2+q^2=0$ on these two  open subset, where $q$ is  non-zero
real valued on the real set $l_{\infty}^{\sigma}$.
 Hence $Z_0$ does not have real
point over
$l_{\infty}$. Thus we have $(Z_0)^{\sigma}=\Phi_0^{-1}(P_0)$.

Finally we show that $\Gamma_0=\mu^{-1}(P_0)$ has no real point. To show this,
we use a coordinate
$(w_0,w_2,w_3)$ (around $P_0$) defined in (\ref{coord}). Then $B$
is  given by  $w_1^2+w_2w_3=0$. Further, it is easy to
see that  the real structure is also given by
$\sigma(w_1,w_2,w_3)=(\ol{w}_1,\ol{w}_3,\ol{w}_2)$.  Now because $\Gamma_0$
is the exceptional curve of a small resolution of an ordinary double
point, $\Gamma_0$ can be canonically identified with the set of lines
contained in the cone
$w_1^2+w_2w_3=0$. 
If $\{(w_1:w_2:w_3)=(a_1:a_2:a_3)\}$ is a real line, we can suppose
$a_1\in\mathbf R$, $a_3=\ol{a}_2$.
It follows that it cannot be 
contained in the cone. This implies that $\Gamma_0$ has no real point.
On the other hand, resolutions of the singularities over $P_{\infty}$ and 
$\ol{P}_{\infty}$ do not yield real points.
Thus we can conclude that $Z$ has no real point.
 \proofend

\vspace{2mm}
Here we summarize necessary conditions for our threefolds to be (birational to) a 
twistor space:

\begin{prop}\label{prop-necessa}
Let $Z$, $B$, $Q$,  $a$ and $b$ be as in Proposition \ref{prop-def-B},
and put $ f(\lambda)=\lambda(\lambda+1)(a\lambda-b)$ as in Proposition
\ref{prop-pos}. 
Then we have: (i) $Q(\lambda,1)^2-f(\lambda)\geq 0$ holds on $\mathbf
R$ and the equality holds iff $\lambda=\lambda_0$, 
(ii) if $f(\lambda)\geq 0$ and $\lambda\neq\lambda_0$ then
$Q(\lambda,1)>\sqrt{f(\lambda)}$. 
\end{prop}

Proposition \ref{prop-def-B} has the following consequence:
\begin{prop}\label{prop-class} Let $g$ be a self-dual metric on
$3\mathbf{CP}^2$ of positive scalar curvature with a non-trivial Killing
field, and assume that 
$g$ is not conformally isometric to LeBrun metric. Then the naturally
induced
$U(1)$-action   on
$3\mathbf{CP}^2$ is uniquely determined up to equivariant diffeomorphisms.
\end{prop}

\noindent Proof. Let $Z$ be the twistor space of $g$. Then $Z$ is as in
Proposition
\ref{prop-def-B}. Let $H_i$ and $C_i$ ($1\leq i\leq 4$) be  as in the
proof of Proposition
\ref{prop-def-B}. Namely, $H_i$ is a real $U(1)$-invariant plane
such that 
$B|_{H_i}$ is a trope  whose reduction is denoted by
$C_i$.  Then $\Phi^{-1}_0(H_i)$ consists of two irreducible components,
both of which are biholomorphic to $H_i$ ($\simeq\mathbf{CP}^2$). Since
$\mu:Z\ra Z_0$ is small,
$\Phi^{-1}(H_i)$ also consists of two irreducible components, which are
denoted by
$D_i$ and $\ol{D}_i$. The natural morphism
$D_i\ra H_i$ is clearly birational.
We now claim that the set $\{D_i,\ol{D}_i\}_{i=1}^4$ of
smooth $U(1)$-invariant divisors on $Z$ is independent of the
choice of small resolutions of $Z_0$.
To see this, recall that if we use  an affine coordinate
$(x_0,x_1,x_2)$ valid on $\{y_3\neq 0\}$
as in the proof of Proposition
\ref{sing. of B}, $Z_0$ is defined by 
$z^2+(x_2+Q(x_0,x_1))^2-x_0x_1(x_0+x_1)(ax_0-bx_1)=0$. So if we put
$\xi=z+i(x_2+Q(x_0,x_1))$ and $\eta=z-i(x_2+Q(x_0,x_1))$,
we get $$Z_0:\hspace{3mm}\xi\eta=x_0x_1(x_0+x_1)(ax_0-bx_1),$$
and the origin corresponds to a compound
$A_3$-singularity $p_{\infty}:=\Phi_0^{-1}(P_{\infty})$. Then the
irreducible components of
$\Phi^{-1}_0(D_i)$  are defined by
$\xi=\ell_i=0$ and $\eta=\ell_i=0$, where $\ell_i$ is one of
$x_0,x_1,x_0+x_1$ and $ax_0-bx_1$. 
 From this, and from the explicit description of small resolutions which
will be explained in \S \ref{ss-special}, we can deduce that 
both of the birational
morphisms $D_i\ra H_i$ and $\ol{D}_i\ra H_i$ are the composition of three blowing-ups,
and that each blowing-up is always performed at just one of the two
$U(1)$-fixed points on the proper transforms of $C_i$ (or
$\Phi^{-1}_0(C_i)$, more precisely). Evidently there are just $2^3=8$
choices  of blowing-ups satisfying this property in all, and the set of
the resulting surfaces is just
$\{D_i,\ol{D}_i\}_{i=1}^4$. Thus $\{D_i,\ol{D}_i\}_{i=1}^4$ is
independent of the choice of a small resolution of $p_{\infty}$.

Since
$D_i+\ol{D}_i$ is a fundamental divisor, and since we  have
$-(1/2)K_Z\cdot L=2$ ($L$ is a twistor line), we have
$D_i\cdot L=\ol{D}_i\cdot L=1$. Hence by a result of Poon
\cite{P92},  $L_i:=D_i\cap \ol{D}_i$ is a twistor line which is
obviously $U(1)$-invariant, and that $L_i$ is contracted to a
point by the twistor fibration $Z\ra 3\mathbf{CP}^2$ which is
$U(1)$-equivariant.
Hence the $U(1)$-action on $3\mathbf{CP}^2$ can be read  from that on
$3\mathbf{CP}^2$. Therefore the conclusion of the proposition follows.
\proofend

\vspace{2mm} The proposition implies that, up to equivariant
diffeomorphisms,  there are just two effective $U(1)$-actions on
$3\mathbf{CP}^2$ which can be the identity component of the isometry
group of a self-dual metric whose scalar curvature is positive. One is
the semi-free $U(1)$-action, which is the identity component of generic
LeBrun metric, and the other is the action  obtained in
Proposition \ref{prop-class}. Of course, there are many other
differentiable $U(1)$-actions: for example, we can get an infinite
number of mutually inequivariant $U(1)$-actions  by first taking an
effective
$U(1)\times U(1)$-action on $3\mathbf{CP}^2$ and then choosing a
$U(1)$-subgroup of $U(1)\times U(1)$.

\begin{prop}\label{prop-pic}
Let $B$ be the quartic surface defined by the equation (\ref{eqn-B}) and suppose 
that $Q$ and $f$ satisfy the assumption in Proposition \ref{prop-realpoint}.
Let $\Phi_0:Z_0\ra\mathbf{CP}^3$ be the double covering branched along $B$, and
$\mu:Z\ra Z_0$ any small resolution (which exists by Propositions \ref{sing. of B}
 and \ref{prop-pos}). Put $\Phi=\mu\cdot\Phi_0$. Then we
have (i) $K_Z\simeq \Phi^*O(-2)$, (ii) the line bundle $(-1/2)K_Z$ is 
uniquely determined.
\end{prop}

\noindent Proof.
Let $K_{Z_0}$ denote the canonical sheaf of $Z_0$. Then we have $K_{Z_0}\simeq
\Phi_0^*(K_{\mathbf{CP}^3}+(1/2)O(B))\simeq\Phi_0^*O(-2)$. Moreover,
since $\mu$ is small, we have $K_Z\simeq \Phi^* K_{Z_0}$. Hence we get
$K_Z\simeq \Phi^*O(-2)$. For (ii) it suffices to show that $H^1(O_Z)=0$.
Since the singularities of $Z_0$ are normal, and since the exceptional curves
of $\mu$ are rational, 
we get by Leray spectral sequence $H^1(O_Z)\simeq H^1(O_{Z_0})$.
Then applying the spectral sequence to $\Phi_0$, and using $\Phi_{0*}O_{Z_0}
\simeq O\oplus O(-2)$ and $R^q\Phi_{0*}O_{Z_0}=0$ for $q\geq 1$, we get
$H^1(O_{Z_{0}})=0$.
\proofend
 
\section{Defining equations of real touching conics}\label{s-detc}
Let $Z$, $\Phi:Z\ra\mathbf{CP}^3$, $B\subset\mathbf{CP}^3$,
$\Phi_0:Z_0\ra\mathbf{CP}^3$ and $\mu:Z\ra Z_0$ have the meaning of the
previous section. Recall that $B$ has just three singular points
which are denoted by $P_0,P_{\infty}$ and $\ol{P}_{\infty}$, where
$P_0$ is a real ordinary double point, and $P_{\infty}$ and
$\ol{P}_{\infty}$ are conjugate pair of elliptic singularities of type
$\tilde{E}_7$. Set $\Gamma_0=\mu^{-1}(P_0)$, which is a real smooth rational
curve. In this section we first study the images of real lines under the map
$\Phi$. (We do not assume that $Z$ is a twistor space.)

\begin{definition} A real irreducible conic $C$ in $\mathbf{CP}^3$ is
called a real touching conic of $B$ if $C\subset B$ or if the
intersection number with
$B$ is at least two for any intersection points.
\end{definition}

\begin{prop}\label{prop-image} Let $L$ be a
real line in $Z$, where line means a smooth rational curve whose normal bundle
in $Z$ is isomorphic to $O(1)^{\oplus 2}$. Then
$\Phi(L)$ is a line in
$\mathbf{CP}^3$ iff
$L\cap \Gamma_0\neq\phi$. Otherwise $\Phi(L)$ is a real touching conic of $B$.
\end{prop}

\noindent Proof. By adjunction formula, we have $-2=K_Z\cdot L+\deg N_{L/Z}=
K_Z\cdot L+2$. Hence we have $(-1/2)K_Z\cdot L=2$. Therefore
 $\Phi(L)$ is a curve whose  degree is at most two,
and
$\Phi(L)$ is a line iff
$\Phi|_L:L\ra\Phi(L)$ is two-to-one. Assume that $\Phi(L)$ is a line,
which is necessarily real. Consider the pencil of planes whose base
locus is $\Phi(L)$.
By Bertini's theorem, general member of this pencil is singular precisely when
$\Phi(L)$ goes through the singular point of $B$. If $P_0\not\in\Phi(L)$,
$\Phi(L)$ is the line joining $P_{\infty}$ and $\ol{P}_{\infty}$.
Namely, $\Phi(L)=l_{\infty}=\{y_0=y_1=0\}$. As we have already seen in the 
proof of Proposition \ref{prop-realpoint}, there is no real point on $\Phi_0^{-1}
(l_{\infty})$ and hence so is  $\Phi^{-1} (l_{\infty})$. 
Therefore $\Phi^{-1} (l_{\infty})$ has no real components.
But because we have chosen real $L$, this is a contradiction.
Hence we have $P_0\in\Phi(L)$.
It follows that $L\cap \Gamma_0\neq\phi$.
Conversely assume that $L$ is a real line intersecting $\Gamma_0$.
Then since there are no real points on $\Gamma_0$
(Proposition \ref{prop-realpoint} (ii)), 
the intersection is not one point. Because
$\Phi(\Gamma_0)=P_0$, this implies that $\Phi$ is not one-to-one on $L$.
Hence $\Phi(L)$ must be a line.

Finally suppose that $\Phi(L)$ is a conic. If $(\Phi(L), B)_P=1$ for
some $P\in \Phi(L)\cap B$, then $P$ is a smooth point of $B$
and  the intersection is transversal. Therefore
$\Phi^{-1}(\Phi(L))$ is locally irreducible near $\Phi^{-1}(P)$.
This contradicts to the fact that $\Phi|_L$ is bijective.
Therefore we have $(\Phi(L), B)_P\geq 2$ for any $P\in \Phi(L)\cap B$,
 which implies that
$\Phi(L)$ is a touching conic of $B$, which is necessarily real.
\hfill $\square$

 \vspace{2mm}
In general, it seems to be impossible to determine the defining equations  of
touching conics of a smooth quartic curve.
Namely, it seems to be  hopeless to write down the equation of touching conics
contained in a general plane.
In the rest of this section we determine defining equations of real touching
conics contained in  real $U(1)$-invariant planes.
By Proposition \ref{prop-def-B} (iii) such a plane is of the form 
$H_{\lambda}=\{y_0=\lambda y_1\}$ or $H_{\infty}=\{y_1=0\}$,
where $\lambda$ is a real number.
As seen in the proof of Proposition \ref{prop-realpoint}, $B_{\lambda}=B\cap
H_{\lambda}$ is a union of two 
$U(1)$-invariant conics, and their
intersection is
$P_{\infty}$ and $\ol{P}_{\infty}$.

Let $C$ be a real touching conic contained in $H_{\lambda}$.
Since  the two components of $B_{\lambda}$ have the same tangent lines at
$P_{\infty}$ and $\ol{P}_{\infty}$,  the local intersection number
$(C, B_{\lambda})_{P_{\infty}}$ is
zero, two, or four, and by reality, the same holds for  $(C,
B_{\lambda})_{\ol{P}_{\infty}}$. Correspondingly, real touching conics can be
classified into the following three types:

\begin{definition}{\em{ Let $C$  be a real touching conic in  $H_{\lambda}$.
Then 
(i)  $C$  is called {\em generic type}\,
if $(C, B_{\lambda})_{P_{\infty}}=(C,B_{\lambda})_{\ol{P}_{\infty}}=0$.
(ii)  $C$  is called {\em special type}
if $(C, B_{\lambda})_{P_{\infty}}=(C,B_{\lambda})_{\ol{P}_{\infty}}=2$.
(iii)  $C$  is called {\em orbit type}
if $(C, B_{\lambda})_{P_{\infty}}=(C,B_{\lambda})_{\ol{P}_{\infty}}=4$.
}}
\end{definition}


If $C$ is a real touching conic of generic type, then $P_{\infty},
\ol{P}_{\infty}\not\in C$ and 
 $C\cap B_{\lambda}$ consists of just four points, all satisfying
$(C,B_{\lambda})_P=2$. 
If $C$ is a real touching conic of special type, then $C$ goes through  $P_{\infty}$
and $
\ol{P}_{\infty}$ but the tangent lines  at
$P_{\infty}$ and $\ol{P}_{\infty}$ are different from those  of $B_{\lambda}$.
Further,  there  are other intersection $P$ and $\ol{P}$ satisfying
$(C,B_{\lambda})_P=(C,B_{\lambda})_{\ol{P}}=2$. 
If $C$ is a real touching conic of orbit type, there are no other intersection
points.  In this case, $C$ is the closure of $\mathbf C^*$-action,
where $\mathbf C^*$-action is the complexification of $U(1)$-action.

First we classify real touching conics of generic type:

\begin{prop}\label{prop-a} (i) If $f(\lambda)>0$ and 
$Q(\lambda,1)^2>f(\lambda)$, there exists a family of real touching conics
of generic type  on
$H_{\lambda}$, parametrized by a circle. Their  defining equations are
explicitly given by
\begin{equation}\label{eqn-a-0}
2(Q^2-f)y_1^2+\sqrt{f}e^{i\theta}y_2^2+2Qy_2y_3+\sqrt{f}e^{-i\theta}y_3^2=0,
\end{equation} where we put  $Q=Q(\lambda,1)$ and $f=f(\lambda)$, and
$\theta\in\mathbf R$. Further, every real touching conic of generic type in
$H_{\lambda}$ is a member of this family. (ii) If $f(\lambda)<0$ or
$Q(\lambda,1)^2=f(\lambda)$, there is no real touching conic of generic
type on
$H_{\lambda}$. (iii) If $Q(\lambda,1)^2>f(\lambda)>0$, the conic
(\ref{eqn-a-0}) has no real point for any
$\theta\in\mathbf R$.
\end{prop}

We note that  $U(1)$ acts transitively (but non-effectively) on the space of these
touching  conics. To prove the proposition, we need the following 

\begin{lemma}\label{lemma-2dr} An equation
$x^4+a_1x^3+a_2x^2+a_3x+a_4=0$ has two double roots iff (i) if $a_1\neq
0$, then $4a_1a_2=a_1^3+8a_3$ and
$a_1^2a_4=a_3^2$ hold, (ii) if $a_1=0$, then  $a_3=0$ and $4a_4=a_2^2$
hold.
\end{lemma}

\noindent Proof. Suppose first that the equation has two double roots
$\alpha$ and $\beta$. Then we have
$x^4+a_1x^3+a_2x^2+a_3x+a_4=(x-\alpha)^2(x-\beta)^2$. Expanding the
right hand side, and comparing the coefficients, we get
$a_1=-2(\alpha+\beta)$, $a_2=\alpha^2+\beta^2+4\alpha\beta$,
$a_3=-2\alpha\beta(\alpha+\beta)$, and $a_4=\alpha^2\beta^2$.  From the
first one, we get $\alpha+\beta=-a_1/2$. If $a_1\neq 0$, we further get
$\alpha\beta=a_3/a_1$. Hence we have
$a_2=(\alpha+\beta)^2+2\alpha\beta=a_1^2/4+2a_3/a_1$ and
$a_4=(\alpha\beta)^2=a_3^2/a_1^2$. From these, we get the equalities of
(i). If
$a_1=0$, it immediately follows from $\alpha+\beta=0$ that $a_3=0$,
$a_2=-2\alpha^2$ and $a_4=\alpha^4$. Hence we get the equalities of
(ii). Thus we obtain the necessity. Conversely, assume $a_1\neq 0$ and
the equalities of (i) holds. Put
$$
\alpha=-\frac{a_1}{4}+\sqrt{\frac{a_1^2}{16}-\frac{a_3}{a_1}},\hspace{5mm}
\beta=-\frac{a_1}{4}-\sqrt{\frac{a_1^2}{16}-\frac{a_3}{a_1}}.
$$ Then it is a straightforward calculation to see that 
$(x-\alpha)^2(x-\beta)^2=x^4+a_1x^3+a_2x^2+a_3x+a_4$ under our assumptions.
 Thus the
given equation has  two double roots $\alpha$ and $\beta$. Finally assume
$a_1=a_3=0$ and  
$4a_4=a_2^2$. Then we have
$x^4+a_1x^3+a_2x^2+a_3x+a_4=x^4+a_2x^2+(a_2^2/4)=(x^2+(a_2/2))^2$. Hence
the equation has two double roots. \proofend

\vspace{2mm}
\noindent Proof of Proposition \ref{prop-a}. Let 
$ay_1^2+by_1y_2+cy_2^2+dy_1y_3+ey_2y_3+hy_3^2=0$ be 
a defining equation of $C$ in $H_{\lambda}$. Then since
$B_{\lambda}\cap\{y_3=0\}=\{\ol{P}_{\infty}\}$ and since $C$ is assumed to be
generic type, all of the touching points are  on  $\{y_3\neq0\}$.
Putting $x_1=y_1/y_3$ and $x_2=y_2/y_3$ on $\{y_3\neq 0\}$ as before, 
$C$ is defined by
\begin{equation}\label{eqn-a}
 ax_1^2+bx_1x_2+cx_2^2+dx_1+ex_2+h=0
\end{equation} 
and $B_{\lambda}$ is defined by (as in the proof of Proposition
\ref{prop-realpoint})
$$\left(x_2+\left(Q-\sqrt{f}\right)x_1^2\right)
\left(x_2+\left(Q+\sqrt{f}\right)x_1^2\right)=0.$$
 Let $g$ denote $g_+:=Q+\sqrt{f}$ or
$g_-:=Q-\sqrt{f}$. Substituting $x_2=-gx_1^2$ into (\ref{eqn-a}), we get
\begin{equation}\label{eqn-a-1} g^2cx_1^4-gbx_1^3+(a-ge)x_1^2+dx_1+h=0.
\end{equation} If $c=0$, (\ref{eqn-a-1}) cannot have two double roots,
so we have $c\neq 0$. Suppose $b\neq 0$. Then by Lemma \ref{lemma-2dr}
(i), (\ref{eqn-a-1}) has two double roots iff 
$$-\frac{4b}{cg}\cdot\frac{a-eg}{cg^2}=-\frac{b^3}{c^3g^3}+\frac{8d}{cg^2}
\hspace{3mm}{\rm{and}}\hspace{2mm}
\left(-\frac{b}{cg}\right)^2\cdot\frac{h}{cg^2}=\frac{d^2}{c^2g^4}
$$ hold. These can be respectively written
\begin{equation}\label{eqn-a-2} 4bc(a-eg)=b^3-8gc^2d,\,\,\,\,\,b^2h=cd^2.
\end{equation} Namely, a conic (\ref{eqn-a}) with $b\neq 0$ is a
touching conic of generic type iff (\ref{eqn-a-2}) is satisfied for both of
$g=g_+$ and $g=g_-$. In this case, simple calculations show that 
$4ac=b^2$ and $4ah=d^2$ and $2ae=bd$. From these we get
$a=b^2/4c, e=2cd/b$ and $h=cd^2/b^2$. Substituting these into (\ref{eqn-a}),
we get $(bx_1+2cx_2+2cd/b)^2=0$. This implies that $C$ is a double line.
 Thus
contradicting our assumption and we get $b=0$. Then by Lemma
\ref{lemma-2dr} (ii), we have $d=0$ and 
\begin{equation}\label{eqn-a-3} 4g^2ch=(a-ge)^2.
\end{equation} If we regard (\ref{eqn-a-3}) as a homogeneous equation
about
$(a:c:e:h)\in\mathbf{CP}^3$, (\ref{eqn-a-3}) is a quadratic cone whose
vertex is
$(a:c:e:h)=(g:0:1:0)$.
We need to get the intersection of these two quadrics. Restricting
(\ref{eqn-a-3}) onto the plane
$a=\kappa e$, we get
\begin{equation}\label{eqn-a-4} 4g^2ch=(\kappa-g)^2e^2,\,\, g=g_{\pm}.
\end{equation} It is readily seen that these two conics (for the case
$g=g_+$ and
$g=g_-$) coincide iff
$\kappa=0$ or $\kappa=(Q^2-f)/Q$. If $\kappa=0$, we have $a=0$, so 
(\ref{eqn-a}) becomes
$cx_2^2+ex_2+h=0$, where the coefficients are subjected to $4ch=e^2$. This
implies that
$C$ is a union of two lines, contradicting  our assumption. Hence we have
$\kappa=(Q^2-f)/Q$.  Then (\ref{eqn-a-4}) becomes 
\begin{equation}\label{eqn-a-5} 4Q^2ch=fe^2.
\end{equation} If $e=0$, it follows that $h=0$, and hence
(\ref{eqn-a}) will again be a union of lines. Therefore we have $e\neq 0$. It is
easily seen that the real structure on the space of coefficients is
given by
$(a:c:e:h)\mapsto (\ol{a}:\ol{h}:\ol{e}:\ol{c})$.   Hence if $f=f(\lambda)<0$,
 (\ref{eqn-a-5})
cannot hold for real $(a:c:e:h)$. Namely, on $H_{\lambda}$, there does not
exist a real touching conic of generic type if
$f(\lambda)<0$. This implies (ii) of the proposition for the case
$f(\lambda)<0$. If $f=f(\lambda)>0$, putting $e=1$ and $h=\ol{c}$ in
(\ref{eqn-a-5}) we get
$4Q^2|c|^2=f$. Hence we can write $$c=\frac{\sqrt{f}}{2Q}e^{i\theta}$$
for some
$\theta\in\mathbf R$. Further we have $a=(Q^2-f)/Q$. Substituting these
into (\ref{eqn-a}), we get (\ref{eqn-a-0}).  Then it is immediate to see
that the determinant of the matrix defining (\ref{eqn-a-0}) is
$-2(Q^2-f)^2$. Therefore the conic (\ref{eqn-a-0}) is irreducible iff
$Q^2-f\neq 0$. Thus we get 
 (i), and also (ii) for the case $Q(\lambda)^2=f(\lambda)$.

Finally we show (iii). Recall that the real structure on $H_{\lambda}$
is given by
$(y_1:y_2:y_3)\mapsto (\ol{y}_1:\ol{y}_3:\ol{y}_2)$. Hence if
$(y_1:y_2:y_3)\in H_{\lambda}$ is a real point, we can suppose $y_1\in\mathbf
R$ and $y_3=\ol{y}_2$. Substituting these into (\ref{eqn-a-0}), we get  
\begin{equation}\label{eqn-a-6}
(Q^2-f)y_1^2+Q|y_2|^2+\sqrt{f}\cdot{\rm{Re}}(e^{i\theta}y_2^2)=0,
\end{equation} where ${\rm{Re}}(z)$ denotes the real part of $z$.
From this it follows that $y_2=0$ implies $y_3=y_1=0$, so we have
$y_2\neq 0$. Then recalling that $Q>\sqrt{f}$ (Proposition \ref{prop-realpoint}
(i)) we have
$Q|y_2|^2>\sqrt{f}|y_2|^2$. Also we have
${\rm{Re}}(e^{i\theta}y_2^2)\leq |y_2|^2$. Therefore we have
\begin{eqnarray*}
(Q^2-f)y_1^2+Q|y_2|^2+\sqrt{f}\cdot{\rm{Re}}(e^{i\theta}y_2^2)& >&
(Q^2-f)y_1^2+\sqrt{f}|y_2|^2-\sqrt{f}|y_2|^2\\ &=& (Q^2-f)y_1^2\,\,\geq 0
\end{eqnarray*} This implies that (\ref{eqn-a-6}) does not hold for any
real $(y_1:y_2:y_3)$ and
$\theta\in\mathbf R$. Thus there is no real point on the conic (\ref{eqn-a-0})
provided $Q^2>f>0$, and we get (iii).
\hfill $\square$

\vspace{2mm}
Next we classify real touching conics of special type:

\begin{prop}\label{prop-b} (i) If $f(\lambda)>0$, there is no real
touching conic of special type on
$H_{\lambda}$. (ii) If $f(\lambda)<0$, there exists a family of 
 real touching conics of special type, parametrized by a circle. Their
defining equations are given by
\begin{equation}\label{eqn-b-0}
\sqrt{Q^2-f}\cdot y_1^2+\sqrt{\frac{\sqrt{Q^2-f}-Q}{2}}\cdot
e^{i\theta}y_1y_2+
\sqrt{\frac{\sqrt{Q^2-f}-Q}{2}}\cdot e^{-i\theta}y_1y_3+y_2y_3=0,
\end{equation} where we put  $Q=Q(\lambda,1)$ and $f=f(\lambda)$,
and $\theta\in \mathbf R$ as before.
Further, every
real touching conic of special type  in $H_{\lambda}$ is a member of 
this family.  (iii) The conic (\ref{eqn-b-0}) has no real point for any
$\theta\in\mathbf R$.
\end{prop}

Note again that $U(1)$ acts transitively on the parameter space of
these touching conics.
Also note that if $f<0$ we have $Q^2-f>0$ and $\sqrt{Q^2-f}-Q>0$. Hence
every square root in the equation make a unique sense (i.e. we always take the positive
root).

\vspace{2mm}
\noindent Proof. Let $C$ be a real touching conic of special type on
$H_{\lambda}$. 
Then since $C$ goes through $P_{\infty}$ and $\ol{P}_{\infty}$, 
the other two touching points
belong to mutually different irreducible component of $B_{\lambda}$.
On the other hand, as shown in the proof of Proposition \ref{prop-realpoint},
each irreducible components of $B_{\lambda}$ is real iff $f(\lambda)\geq 0$.
Therefore, on $H_{\lambda}$, there does not exist real touching conic of special type 
if $f(\lambda)>0$. Thus we get (i). So in the sequel we suppose $f(\lambda)<0$.

It is again readily seen that touching points are not on the line
$\{y_3=0\}$.  So we still use $(x_1,x_2)$ as a non-homogeneous coordinate
 on
$\mathbf C^2=\{y_3\neq 0\}\subset H_{\lambda}$. Then because
$C$ contains $P_{\infty}$ and
$\ol{P}_{\infty}$, an equation of a touching conic $C$ of special
type is of the form
\begin{equation}\label{eqn-b} ax_1^2+bx_1x_2+dx_1+ex_2=0.
\end{equation} Substituting $x_2=-gx_1^2$ into (\ref{eqn-b}), we get
\begin{equation}\label{eqn-b-1}
x_1\cdot\left(gbx_1^2+(ge-a)x_1-d\right)=0.
\end{equation} If $d=0$, $x_1=0$ is a double root of (\ref{eqn-b-1}).
Then the tangent line of $C$ at $P_{\infty}=(0,0)$ becomes $x_2=0$,  as
is obvious from (\ref{eqn-b}). This implies that $C$ is a touching conic of
orbit type, contradicting our assumption.Therefore, we have $d\neq 0$. (\ref{eqn-b-1}) has a
double root other than $x_1=0$ iff
\begin{equation}\label{eqn-b-2} (ge-a)^2+4gbd=0.\end{equation} Namely
(\ref{eqn-b}) is a touching conic of special type iff (\ref{eqn-b-2}) is
satisfied for both of $g=g_+$ and $g=g_-$. (Note that $g_+\neq g_-$.) If we
regard (\ref{eqn-b-2}) as a  homogeneous equation of
$(a:b:d:e)\in\mathbf{CP}^3$, (\ref{eqn-b-2}) is a quadratic cone whose
vertex is 
$(a:c:d:e)=(g:0:0:1)$. That is, the parameter space of  touching conics
of special type is the intersection of two quadratic cones in
$\mathbf{CP}^3$. Restricting (\ref{eqn-b-2}) onto the plane $a=\kappa
e$, we get
$$ (g-\kappa)^2e^2+4gbd=0,\,\, g=g_{\pm}.$$ It is readily seen that these
two conics coincide iff $\kappa=\pm\sqrt{g_+g_-} =\pm\sqrt{Q^2-f}$. Therefore $C$ is a
touching conic of special type iff
\begin{equation}\label{eqn-b-3} a=\sqrt{Q^2-f}\cdot e,\,\,\,
\left(Q-\sqrt{Q^2-f}\right)e^2+2bd=0
\end{equation} or
\begin{equation}\label{eqn-b-4} a=-\sqrt{Q^2-f}\cdot e,\,\,\,
\left(Q+\sqrt{Q^2-f}\right)e^2+2bd=0
\end{equation} hold. 
 It is easily seen that the real structure on the space of coefficients
is given by
$(a:b:d:e)\mapsto (\ol{a}:\ol{d}:\ol{b}:\ol{e})$.  

Since we have assumed  $f<0$, we have $Q+\sqrt{Q^2-f}>0$ 
and there is no real conic satisfying (\ref{eqn-b-4}). On the other
hand, we have
$Q-\sqrt{Q^2-f}<0$. Hence by (\ref{eqn-b-3}) we have
$$ 2|b|^2=\left(\sqrt{Q^2-f}-Q\right)e^2,
$$ where $b\in \mathbf C$ and $e\in\mathbf R$. If $e=0$, then $b=a=0$,
contradicting  the assumption that $C$ is a conic. Hence $e\neq 0$, and
we may put $e=1$. Then we can write
$$b=\sqrt{\frac{\sqrt{Q^2-f}-Q}{2}}\cdot e^{i\theta}$$ for some
$\theta\in\mathbf R$. Also we have $d=\ol{b}$. Thus we obtain
(\ref{eqn-b-0}) of the proposition.

Finally we show (iii). If $(y_1:y_2:y_3)$ is a real point of
$H_{\lambda}$, we can suppose $y_1\in\mathbf R$ and $y_3=\ol{y}_2$.
 Substituting
these into (\ref{eqn-b-0}), we get
\begin{equation}\label{eqn-b-5}
\sqrt{Q^2-f}\cdot y_1^2+\sqrt{2\left(\sqrt{Q^2-f}-Q\right)}\cdot
y_1\cdot{\rm{Re}}(e^{i\theta}y_2)+|y_2|^2=0.
\end{equation} If $y_1=0$, it follows $y_2=y_3=0$. Hence $y_1\neq 0$
 and we can suppose $y_1>0$. Then
 we have $y_1{\rm{Re}}(e^{i\theta}y_2)\geq -y_1|y_2|$. Hence we
have
\begin{eqnarray*}
\lefteqn{
\sqrt{Q^2-f}\cdot y_1^2+\sqrt{2\left(\sqrt{Q^2-f}-Q\right)}\cdot
y_1{\rm{Re}}(e^{i\theta}y_2)+|y_2|^2}\\ &\geq&\sqrt{Q^2-f}\cdot
y_1^2-\sqrt{2\left(\sqrt{Q^2-f}-Q\right)}\cdot y_1|y_2|+|y_2|^2\\
&=&\left(|y_2|-\sqrt{\frac{\sqrt{Q^2-f}-Q}{2}}\cdot y_1\right)^2+
\frac{\sqrt{Q^2-f}+Q}{2}\cdot y_1^2
\end{eqnarray*} Because $y_1\neq 0$ and $f<0$, we have $(\sqrt{Q^2-f}+Q)y_1^2>0$.
Therefore, the left hand side of (\ref{eqn-b-5}) is strictly positive. Thus 
(\ref{eqn-b-5}) does not hold for any real $(y_1:y_2:y_3)\in H_{\lambda}$   and any
$\theta\in\mathbf R$. Therefore the conic (\ref{eqn-b-0}) has no real point for any
$\theta\in\mathbf R$.

It is immediate to see that the determinant of the matrix defining
(\ref{eqn-b-0}) is $-(Q+\sqrt{Q^2-f})/8$ and this is negative if $f<0$.
Hence the conic (\ref{eqn-b-0}) is irreducible
\hfill$\square$

\vspace{3mm} The case of orbit type is straightforward and need no assumption
on the sign of $f(\lambda)$:
\begin{prop}\label{prop-c}  There exists a family of real touching
conics of orbit type, parametrized by non-zero real numbers. Their 
defining equations are 
\begin{equation}\label{eqn-c} y_2y_3=\alpha y_1^2,\,\,\alpha\in \mathbf
R^{\times}.
\end{equation} Further, every real touching conic of orbit type in
$H_{\lambda}$ is contained in this family. 
\end{prop}

Note that by Lemma \ref{lemma-element}, the conic (\ref{eqn-c})
has no real point iff $\alpha<0$.

Combining Propositions \ref{prop-a}, \ref{prop-b} and \ref{prop-c}, we get the
following

\begin{prop}\label{prop-type} Let
$\{S_{\lambda}:=\Phi^{-1}(H_{\lambda})\set H_{\lambda}\in\langle
l_{\infty}\rangle^{\sigma}\}$ be the real members of the pencil of
$U(1)$-invariant  divisors on $Z$. Then  (i) if
$f(\lambda)>0$,  the images of real lines in $S_{\lambda}$ are real
touching conics of generic type or orbit type, (ii) if $f(\lambda)<0$, 
the images of real lines in
$S_{\lambda}$ are real touching conics of special type or orbit type.
\end{prop}

The ambiguity of the types in this proposition will be
  excluded in Section
\ref{s-nb} (Theorem \ref{thm-type}), for the images of twistor lines.

\section{The inverse images of real touching conics}\label{s-inv} According
to the previous section, real touching conics in fixed $H_{\lambda}$ form
 families parametrized by a circle for generic and special types, or
$\mathbf R^{\times}$ for orbit type. In this section we study the
inverse images of these touching conics in $Z$, which are
 candidates of  twistor lines. We begin with the following

\begin{prop}\label{prop-inv-pl} If $\lambda\in\mathbf R$ and if 
$f(\lambda)\neq 0$, then
$S_{\lambda}:=\Phi^{-1}(H_{\lambda})$ is a smooth rational surface with
$c_1^2=2$, and is a real $U(1)$-invariant member of $|(-1/2)K_Z|$.
\end{prop}
\noindent Proof.
 Smoothness of $S_{\lambda}$ can be checked by giving all small
resolutions of the singular threefold
$Z_0$ (which is the double cover branched along $B$). Indeed,
if $f(\lambda)\neq 0$, $P_{\infty}$ and 
$\ol{P}_{\infty}$ are
$A_3$-singularities of the curve $B_{\lambda}=H_{\lambda}\cap B$, so that the surface
$\Phi^{-1}_0(H_{\lambda})$ has $A_3$-singularities over there, which are
minimally resolved through small resolutions of
$Z_0$.
If $\lambda=\lambda_0$, $P_0=(\lambda_0:1:0:0)$ is a node of $B_{\lambda_0}$,
but this is also resolved by small resolutions of the corresponding ordinary double 
points of $Z_0$.
See \S\ref{ss-special} and \S\ref{s-line}, where  all resolutions of
singularities of $Z_0$ are concretely given. In order to see the structure of
$S_{\lambda}$, note that  $S_{\lambda}$
 is a smooth deformation of a surface $\Phi^{-1}(H)$, where $H$ is a
general plane in $\mathbf{CP}^3$. It is easy to that $\Phi^{-1}(H)$ is a
smooth rational surface with $c_1^2=2$, because  $B\cap H$ is 
a smooth quartic curve for general $H$. From
this it follows the same properties for
$S_{\lambda}$,since they can be connected by smooth deformations.
 Then by Proposition \ref{prop-pic} we have
$S_{\lambda}\in|\Phi^*O(1)|=|(-1/2)K_Z|$. The reality and the $U(1)$-invariance of
$S_{\lambda}$ is obvious from our choice of $H_{\lambda}$.
\proofend

\vspace{2mm}
Next we investigate the inverse images of touching conics of generic type.

\begin{prop}\label{prop-inv-a} Suppose $f(\lambda)>0$ and $\lambda\neq\lambda_0$
(namely $Q^2-f>0$), and let
$C_{\theta}\subset H_{\lambda}$ be a real touching conic of generic type
defined by the equation (\ref{eqn-a-0}).  Then the following (i)--(iv)
hold:  (i) 
$\Phi^{-1}(C_{\theta})$ has just two irreducible components, both of
which are smooth rational curves that are mapped biholomorphically onto
$C_{\theta}$, (ii) each irreducible component of $\Phi^{-1}(C_{\theta})$
has a trivial normal bundle in $S_{\lambda}$,   (iii)  these two
irreducible components of
$\Phi^{-1}(C_{\theta})$ belong to mutually different pencils on
$S_{\lambda}$, (iv) each  irreducible component of 
$\Phi^{-1}(C_{\theta})$ is  real.
\end{prop}

\noindent Proof. Since $C_{\theta}$  and the branch quartic 
$B_{\lambda}$  have the same
tangent line at any intersection points, it is obvious that 
$\Phi_0^{-1}(C_{\theta})$  splits into two irreducible  components
$L_1$ and $L_2$ which are
mapped biholomorphically onto $C_{\theta}$.  Thus we get (i).
 For (ii) first note that
$\Phi^{-1}(C_{\theta})=L_1+L_2$ belongs to $|-2K|$ of
$S_{\lambda}$ since we have $\Phi^*O_{H_{\lambda}}(1)\simeq -K$.
 Hence we have
$(-2K)^2=(L_1+L_2)^2=L_1^2+L_2^2+2L_1L_2$ on $S_{\lambda}$.
 On the other hand, we
have
$4c_1^2=8$  by Proposition
\ref{prop-inv-pl}. Hence we get $L_1^2+L_2^2+2L_1L_2=8$. Further, since $L_1$
and $L_2$ intersect transversally at four points (over the
touching points of
$C_{\theta}$  with $B_{\lambda}$), we have $L_1L_2=4$. Therefore
we get
$L_1^2+L_2^2=0$. Moreover, by (\ref{eqn-a-0}), $C_{\theta}$ actually moves in
a holomorphic family of curves on
$H_{\lambda}$. Hence we have $L_1^2\geq 0$ and $L_2^2\geq 0$. Therefore we get
$L_1^2=L_2^2=0$. Namely we have (ii). (iii) immediately follows from
(ii), since we have $L_1L_2=4$ on $S_{\lambda}$. (iv) is harder than  one may think
at first glance,  since there is no real point on
$C_{\theta}$. First we note that it suffices to prove the claim for
$C_0$ (= the curve obtained by setting $\theta=0$ for $C_{\theta}$), since
$U(1)$ acts transitively on the parameter space of real  touching conics
of generic type (see Proposition \ref{prop-a}).
 The idea of our proof  of the reality is as follows:   the map
$\Phi^{-1}(C_0)\ra C_0$ is  finite, two sheeted covering whose branch
consists of four points.  We choose a real simple closed curve $\mathcal
C$ in $C_0$ containing all of these branch points,  in such a way that
over $\mathcal C$  we can
distinguish two sheets, so that 
 we can explicitly see the reality of each irreducible components. To
this end, we still use $(y_1:y_2:y_3)$ as a homogeneous coordinate on $H_{\lambda}$
and set $U:=\{y_1\neq 0\}$ ($\simeq \mathbf{C}^2$), which is clearly a
real subset of
$H_{\lambda}$, and use $(v_2,v_3)=(y_2/y_1,y_3/y_1)$ as an affine coordinate on
$U$. Then $Z_0|_{U}=\Phi^{-1}_0(U)$ is
defined by the equation
\begin{equation}\label{eqn-Z-U} z^2+\left(v_2v_3+Q\right)^2-f=0,
\end{equation} where  $z$ is a fiber coordinate of $ O(2)$, and
$Q=Q(\lambda,1)$, $f=f(\lambda)$ as in the previous section.  The real
structure is given by 
$(v_2,v_3, z)\mapsto (\ol{v}_3,\ol{v}_2,\ol{z})$. 
  Then on $U$, our equation (\ref{eqn-a-0}) of $C_0$ becomes 
$2(Q^2-f)+\sqrt{f}v_2^2+2Qv_2v_3+\sqrt{f}v_3^2=0$.  Now we introduce a
new coordinate $(u,v):=(v_2+v_3,v_2-v_3)$ which is valid on $U$. Then our real
structure is given by $(u,v)\mapsto (\ol{u},-\ol{v})$, and the defining
equation of $C_0$ becomes
$4(Q^2-f)+\sqrt{f}(u^2+v^2)+Q(u^2-v^2)=0$. From this we immediately have 
\begin{equation}\label{eqn-uv-1}
C_0:\hspace{3mm}v^2=4\left(Q+\sqrt{f}\right)+\frac{Q+\sqrt{f}}{Q-\sqrt{f}}u^2.
\end{equation} We put $V:=\{(u,v)\in  U\set u\in i\mathbf{R},
v\in\mathbf R\}$ which  is clearly a real subset of $U$. Then $\mathcal
C:=V\cap C_0$ is a real simple closed curve (an ellipse) in $V\simeq\mathbf R^2$.
 Indeed,
putting $u=iw$ ($w\in\mathbf R$), we get from (\ref{eqn-uv-1})  
$$C_0:\hspace{3mm}v^2+\frac{Q+\sqrt{f}}{Q-\sqrt{f}}w^2=4\left(Q+\sqrt{f}\right).$$
(Note that  $Q-\sqrt{f}>0$ by our assumption $f>0$ and Proposition
\ref{prop-necessa}.)
  Substituting
$v_2v_3=(u^2-v^2)/4$ into (\ref{eqn-Z-U}), and then using (\ref{eqn-uv-1}), 
we get
\begin{equation}\label{eqn-inv-C}
\Phi^{-1}(C_0):\hspace{3mm}4\left(Q-\sqrt{f}\right)^2z^2
+fu^2\left(u^2+4\left(Q-\sqrt{f}\right)\right)=0,
\end{equation} or,  using $w$ above, 
\begin{equation}\label{eqn-inv-C2}
\Phi^{-1}(C_0):\hspace{3mm}4\left(Q-\sqrt{f}\right)^2z^2
=fw^2\left(4\left(Q-\sqrt{f}\right)-w^2\right).
\end{equation}
 Here note that in (\ref{eqn-uv-1}) $u$ can be used as a coordinate on
$C_0$, only outside the two branch points of $C_{\lambda}\cap U\ra\mathbf C$
defined by
$(u,v)\mapsto u$. In a neighborhood of these branch points,  we
have to use $v$ instead of $u$ as a local coordinate on $C_0$. Then we can
see that the inverse image of a neighborhood of
$u=\pm 2i(Q-\sqrt{f})^{\frac{1}{2}}$ \,(i.e. the branch points) also
splits into two irreducible components, which is of course as expected.
From (\ref{eqn-inv-C}), we easily deduce that the branch points of
$\Phi^{-1}(C_0)\ra C_0$ are $(u,v)=(0,\pm2(Q+\sqrt{f})^{\frac{1}{2}})$
and
$(u,v)=(\pm 2i(Q-\sqrt{f})^{\frac{1}{2}},0)$. All of these four points
clearly lie on
$\mathcal C$, and
$\mathcal C$ is divided into four segments. It immediately follows from
(\ref{eqn-inv-C2}) that 
 $z$ always takes  real value over $\mathcal C$. Moreover, it is clear
that  the sign of $z$ is constant on each of the four segments in
$\mathcal C$, and that the sign changes when passing though the branch
points. On the other hand, since the real structure is given by
$(w,v)\mapsto (-w, -v)$ on $V$, the real structure on
$\mathcal C$ sends each segment to another segment which is not adjacent
to the original one. From these, and because the real structure on
$\Phi^{-1}(\mathcal C)$ is given by $(w,z)\mapsto (-w,\ol{z})=(-w,z)$,
it follows that each of the two irreducible component of
$\Phi^{-1}(\mathcal C)$ is real. Hence the same is true for $\Phi^{-1}(C_0)$.
Thus we have proved (iv) of the proposition.
\hfill
 \proofend

\vspace{2mm}
We have similar statements for touching conics of special type:

\begin{prop}\label{prop-inv-b} Assume $f(\lambda)<0$ and let
$C_{\theta}\subset H_{\lambda}$ be a real touching conic of special type
given by the equation (\ref{eqn-b-0}). Then  (i)--(iv) of Proposition
\ref{prop-inv-a} hold if we replace $\Phi^{-1}(C_{\theta})$ by
$\Phi^{-1}(C_{\theta})-\Gamma-\ol{\Gamma}$, where  we set
$\Gamma:=\Phi^{-1}(P_{\infty})$ and
$\ol{\Gamma}:=\Phi^{-1}(\ol{P}_{\infty})$.
\end{prop}

\noindent Proof. (i) can be proved in the same way as in Proposition
\ref{prop-inv-a}. (But in this case, any small resolution $Z\ra Z_0$
gives the normalization of 
$\Phi^{-1}_0(C_0)$ over $P_{\infty}$ and $\ol{P}_{\infty}$, as will be
mentioned below.) For (ii) first note that we have
$\Phi^{-1}(C_{\theta})=L_1+L_2+\Gamma+\ol{\Gamma}\in |-2K|$ this time,
where
$L_1$ and
$L_2$ are irreducible components of
$\Phi^{-1}(C_{\theta})-\Gamma-\ol{\Gamma}$. 
$\Gamma$ and $\ol{\Gamma}$ are chains of  three $(-2)$-curves on
$S_{\lambda}=\Phi^{-1} (H_{\lambda})$, since they are exceptional curves
of the minimal resolution of $A_3$-singularities of surface.
We write $\Gamma=\Gamma_1+\Gamma_2+\Gamma_3$, where $\Gamma_i$'s are 
smooth rational curves satisfying $\Gamma_1\Gamma_2=\Gamma_2\Gamma_3=1$ and
$\Gamma_1\Gamma_3=0$ on $S_{\lambda}$.
We then have 
$\Gamma^2=\ol{\Gamma}^2=-2$. Furthermore, as we will see in Lemma \ref{lemma-tlr},
 we have $L_1\Gamma=L_1\ol{\Gamma}=L_2\Gamma=L_2\ol{\Gamma}=1$. Therefore
again by Proposition \ref{prop-inv-pl}, we get
$8=(-2K)^2= (L_1+L_2+\Gamma+\ol{\Gamma})^2=L_1^2+L_2^2+2L_1L_2+4$.  But
because
$L_1$ and $L_2$ intersect transversally at two points 
(over the touching points of $C_{\theta}$ and $B_{\lambda}$), and
because $L_1$ and $L_2$ do not intersect on $\Gamma\cup\ol{\Gamma}$,
we have
$L_1L_2=2$. Therefore we have $L_1^2+L_2^2=0$. Hence by the same reason in 
the proof of the previous proposition, we again have
$L_1^2=L_2^2=0$. This implies (ii). (iii) follows from (ii),
because we have $L_1L_2=2$ as is already seen. (iv) can be proved by the
same idea as in the previous proposition: first we may assume
$\theta=0$. Then by (\ref{eqn-b-0}) the equation of
$C_0$ on $U=\{y_1\neq 0\}=\{(v_2,v_3)\}$ is given by
$$
\sqrt{Q^2-f} +\sqrt{\frac{\sqrt{Q^2-f}-Q}{2}}\cdot
(v_2+v_3)+v_2v_3=0.
$$   If we
use another coordinate $(u,v)$ defined in the proof of the previous
proposition, this can be written as
\begin{equation}\label{eqn-circl1}
 C_0:\,\,v^2=u^2+2\sqrt{2}\sqrt{\sqrt{Q^2-f}-Q}\cdot
u+4\sqrt{Q^2-f}.
\end{equation}
 Next we introduce a new variable $w$ by setting 
$u=-\sqrt{2}(\sqrt{Q^2-f}-Q)^{\frac{1}{2}}+iw$. Then the equation becomes
\begin{equation}\label{eqn-bfdj}
C_0:\hspace{3mm}v^2+w^2=2\left(\sqrt{Q^2-f}+Q\right).
\end{equation}
Put $\mathcal C:=C_0\cap\mathbf R^2$, where $\mathbf
R^2=\{(w,v)\set w\in\mathbf R,v\in\mathbf R\}$. Then since 
$\sqrt{Q^2-f}+Q>0$, 
$\mathcal C$ is a real circle in $\{(w,v)\in\mathbf R^2\}$. 

By (\ref{eqn-Z-U}) we have 
\begin{equation}\label{eqn-abrcbr}
\Phi_0^{-1}(U):\hspace{2mm}
z^2=f-\left(\frac{u^2-v^2}{4}+Q\right)^2.
\end{equation}
On the other hand,  by (\ref{eqn-circl1}), we have 
$$C_0: \hspace{2mm}u^2-v^2=-2\sqrt{2}\sqrt{\sqrt{Q^2-f}-Q}\cdot
u-4\sqrt{Q^2-f}
$$ on $C_0$. Substituting this into (\ref{eqn-abrcbr}), we get
$$
\Phi_0^{-1}(C_0):\,\, z^2=f+\frac{\sqrt{Q^2-f}-Q}{2}w^2.
$$
Then by using (\ref{eqn-bfdj}), we get
\begin{equation}\label{eqn-circdl1}
\Phi_0^{-1}(C_0):\,\, z^2=-\frac{\sqrt{Q^2-f}-Q}{2}v^2.
\end{equation}
Therefore, $z$ is pure imaginary over $\mathcal C$,
so that we can distinguish  two sheets by looking the sign of $z/i$.
By (\ref{eqn-circdl1}) and (\ref{eqn-bfdj}), the branch points of
$\Phi^{-1}(C_0)\ra C_0$ is the two points
$(w,v)=(\pm (2(\sqrt{Q^2-f}+Q))^{\frac{1}{2}},0)$
 which lie on $\mathcal C$. The real structure is given by
$(w,v)\mapsto (-\ol{w},-\ol{v})$ and this is equal to $(-w,-v)$ on
$\mathcal C$. Thus the real structure on $\mathcal C$ interchanges the
two segments separated by the two branch points.     
Moreover, the real structure on the fiber coordinate is given by $z\mapsto \ol{z}$. 
Therefore, it changes the sign of $z/i$ over $\mathcal C$.
This implies that each component of $\Phi^{-1}(\mathcal C)$
is real. Therefore that of $\Phi^{-1}(C_0)$ is also real.
\proofend

\vspace{3mm} The situation slightly changes for touching conics of orbit
type:

\begin{prop}\label{prop-inv-c} Suppose $f(\lambda)\neq 0$ and let
$C_{\alpha}\subset H_{\lambda}$ be a real touching conic of orbit type
given by the equation (\ref{eqn-c}).  Then we have: (i) $C_{\alpha}$ is
contained in $B$ iff $\alpha=-Q\pm \sqrt{f}$. (ii)
$\Phi^{-1}(C_{\alpha})-\Gamma-\ol{\Gamma}$ has just two irreducible
components, both of which are smooth rational curves that are mapped
biholomorphically onto
$C_{\alpha}$. (Here $\Gamma$ and $\ol{\Gamma}$ are as in Proposition
\ref{prop-inv-b}.) (iii) Each irreducible component of
$\Phi^{-1}(C_{\alpha})-\Gamma-\ol{\Gamma}$ has a trivial normal bundle in
$S_{\lambda}=\Phi^{-1}(H_{\lambda})$.  (iv) The two irreducible components
of
$\Phi^{-1}(C_{\alpha})-\Gamma-\ol{\Gamma}$ belong to one and  the same
pencil on
$S_{\lambda}$.  (v) Each irreducible component  of
$\Phi^{-1}(C_{\alpha})-\Gamma-\ol{\Gamma}$ is real iff $f>0$ and 
$-Q-\sqrt{f}\leq\alpha\leq-Q+\sqrt{f}$ are satisfied.
(vi) There is no real point on $C_{\alpha}$ if $f$ and $\alpha$ satisfies
the inequalities of (v).
\end{prop}

Note that $f>0$ implies $Q\geq\sqrt{f}$ by Proposition \ref{prop-realpoint}
(ii). 

\vspace{2mm}
\noindent Proof. (i) Substituting $y_2y_3=\alpha y_1^2$ into the defining
equation of $B_{\lambda}$, we get
$\left((\alpha+Q)^2-f\right)y_1^4=0$. Thus if
$C_{\alpha}$ is contained in $ B$ iff 
$(\alpha+Q)^2=f$, which implies $\alpha=-Q\pm \sqrt{f}$, as desired. (ii)
can be seen in the same way as in (i) of Proposition \ref{prop-inv-a}.
(This time, any small  resolution $Z\ra Z_0$ gives the
normalization of 
$\Phi^{-1}_0(C_0)$.)  Next we prove (iii). Let
$\Gamma=\Gamma_1+\Gamma_2+\Gamma_3$, $\ol{\Gamma}=
\ol{\Gamma}_1+\ol{\Gamma}_2+\ol{\Gamma}_3$, $L_1$ and
$L_2$ have the same meaning as in the proof of the last proposition. Then (as
will be shown in Lemma \ref{lemma-intersection3} by explicit calculations)
we have
$L_1L_2=\Gamma\ol{\Gamma}=L_1\Gamma_1=L_1\ol{\Gamma}_1=L_2\Gamma_1=L_2\ol{\Gamma}_1=
L_1\Gamma_3=L_1\ol{\Gamma}_3=L_2\Gamma_3=L_2\ol{\Gamma}_3=0$
 on $S_{\lambda}$, while
$L_1\Gamma_2=L_1\ol{\Gamma}_2=L_2\Gamma_2=L_2\ol{\Gamma}_2=1$
(also on $S_{\lambda}$). In
particular, we have
$L_1\Gamma=L_1\ol{\Gamma}=L_2\Gamma=L_2\ol{\Gamma}=1$.   On the other
hand, we still have $8=(-2K)^2= (L_1+L_2+\Gamma+\ol{\Gamma})^2$
and $L_1^2\geq 0$ and $L_2^2\geq 0$.
Combining these, we get
$L_1^2=L_2^2=0$ on $S_{\lambda}$. Thus we have  (iii). (iv) easily
follows if we consider  the linear systems $|L_1|$ and $|L_2|$, and if we
note that $L_1 L_2=0$. Next we show (v). Substituting $v_2v_3=\alpha$ into
(\ref{eqn-Z-U}), we get
$ z^2+(\alpha+Q)^2-f=0.
$ From this, the equations of irreducible components of 
$\Phi_0^{-1}(C_{\alpha})$ can be  calculated to be
\begin{equation}\label{eqn-inv-gen} z=\pm \sqrt{f-(\alpha+Q)^2}
\end{equation} 
Recalling
that the real structure is given by $z\mapsto \ol{z}$,
these curves
are real iff $f-(\alpha+Q)^2\geq 0$.  In particular, $f\geq 0$ follows.
Then we have
$-Q-\sqrt{f}\leq\alpha\leq -Q+\sqrt{f}$, and we get (v). 
Finally, (vi) immediately follows from Lemma \ref{lemma-element}.
\proofend

\vspace{2mm}
Proposition \ref{prop-inv-c} implies that  not all  touching conics of
orbit type can be the image of a twistor line:   $f(\lambda)>0$
is needed, and further, 
$-Q-\sqrt{f}\leq\alpha\leq-Q+\sqrt{f}$ must be satisfied.
But once we know that one of the two irreducible components is a twistor
line, it follows that  the other component is also a twistor line,  since
by (iv) these two components can be connected by deformation in
$S_{\lambda}$ (hence also in
$Z$) preserving the real structures.
 This is not true for touching conics of
generic type and special type, because the two irreducible components
of $\Phi^{-1}(C_{\theta})$ or $\Phi^{-1}(C_{\theta})-\Gamma-\ol{\Gamma}$
 intersect as
we have already seen in the proofs of Propositions \ref{prop-inv-a} and
\ref{prop-inv-b}.


\section{The normal bundles of the inverse images of real touching conics}
\label{s-nb}
 In
this section we calculate the normal bundle of $L$ in $Z$, where $L$ is
a real irreducible component of the inverse images of the real
touching conics which are determined in Section \ref{s-detc}.
Roughly
speaking, our result states that the normal bundle is isomorphic to
$ O(1)^{\oplus 2}$ for general $L$, but sometimes degenerates into $
O\oplus O(2)$. We can precisely detect which $L$ has such a degenerate
normal bundle, in terms of $Q$ and $f$ which appear in the
defining equation of
$B$.

\subsection{Preliminary lemma and notations} In order to make a distinction between
$O(1)^{\oplus 2}$ and $O\oplus O(2)$,
 we use the following elementary criterion:

\begin{lemma}\label{lemma-nb} Let $E\ra\mathbf{CP}^1$ be a holomorphic
line bundle of rank two, and assume that 
$E$ is isomorphic to $ O(1)^{\oplus 2}$ or $O\oplus O(2)$. Let $s$ and
$t$ be global sections of $E$ which are  linearly independent as global
sections. Then $E\simeq O\oplus O(2)$ iff there are constants
$a,b\in\mathbf C$ such that 
$as+bt$ has two zeros. \end{lemma}

\noindent Proof. It is immediate to see that non-zero sections of $O(1)^{\oplus 2}$ 
have at most one zero. So sufficiency follows. Conversely let
$s=(s_1,s_2)$ and
$t=(t_1,t_2)$ be any linearly independent sections of $O\oplus O(2)$,
where
$s_1,t_1\in\Gamma(O)=\mathbf C$ and
$s_2,t_2\in\Gamma (O(2))$. Then take $a,b\in\mathbf C$ such that
$as_1+bt_1=0$. Then  $as+bt$ can be regarded as a non-zero section of $O(2)$ so
that it has two zeros.
\proofend

\vspace{2mm} Next we introduce some notations.
As in the previous sections, $\lambda\in\mathbf R$ denotes a parameter on 
the space of real $U(1)$-invariant planes.
In other words, $\lambda$ is a parameter on the real locus $l_0^{\sigma}$
of the real line $l_0:=\{y_2=y_3=0\}$.
 The
function $f(\lambda)=\lambda(\lambda+1)(a\lambda-b)$, ($a,b>0$) defines
four open intervals in the circle $l_0^{\sigma}$:
$$I_1=(-\infty,-1),\, I_2=(-1,0),\, I_3=(0,b/a)\hspace{2mm}
{\rm{and}}\hspace{2mm} I_4=(b/a,+\infty). $$ Namely,
$I_1\cup I_3=\{\lambda\in\mathbf R\set f(\lambda)<0\}$ and $I_2\cup
I_4=\{\lambda\in\mathbf R\set f(\lambda)>0\}$. By Proposition
\ref{prop-pos}, the equation $Q(\lambda,1)^2-f(\lambda)=0$ has a
unique real solution
$\lambda=\lambda_0$ which is necessarily a   double root.  Since
we have
$f(\lambda_0)=Q(\lambda_0,1)^2>0$, $\lambda_0\in I_2\cup I_4$. 
By a possible application of a projective transformation with respect to
$y_0$ and
$y_1$, we may suppose that $\lambda_0\in I_4$.  Then we set
$I_{4}^-=(b/a,\lambda_0)$ and $I_{4}^+=(\lambda_0, +\infty)$. 

Next suppose $\lambda\in I_2\cup I_4^-\cup I_4^+$, and let
$$\mathcal C_{\lambda}^{\,\rm{gen}}=\{C_{\theta}\subset H_{\lambda}
\set C_{\theta} {\mbox{ is defined by (\ref{eqn-a-0})}}\}$$ be
the set of real touching conics of generic type on $H_{\lambda}$. Note that if
$\lambda=\lambda_0$ or if $\lambda\in I_1\cup I_3$ there is no real
touching conic of generic type on
$H_{\lambda}$ by Proposition \ref{prop-a} (ii).
Similarly, for $\lambda\in I_1\cup I_3$, 
let $$\mathcal C_{\lambda}^{\,\rm{sp}}=\{C_{\theta}\subset
H_{\lambda}\set C_{\theta} {\mbox{ is defined by (\ref{eqn-b-0})}}\}\hspace{2mm} 
$$ be
the set of real touching conics of special type  on
$H_{\lambda}$.
Note  that if $\lambda\in I_2\cup I_4$ there is no real touching
conic of generic type by Proposition \ref{prop-b} (i). Finally  for
$\lambda\in  I_2\cup I_4$, let
$$\mathcal C_{\lambda}^{\,\rm{orb}}=\left\{C_{\alpha}\subset H_{\lambda}
\set -Q-\sqrt{f}\leq\alpha\leq -Q+\sqrt{f},\hspace{1.5mm}
 C_{\theta} {\mbox{ is defined by (\ref{eqn-c})}}\right\}$$
be
the set of real touching conics of  orbit type on
$H_{\lambda}$.
Note that the restriction on $\alpha$ implies that the two irreducible
components of the inverse image are respectively real (Proposition
\ref{prop-inv-c} (v)). 
$\mathcal C_{\lambda}^{\,\rm{gen}}$ and $\mathcal
C_{\lambda}^{\,\rm{sp}}$ are parametrized by a circle on which $U(1)$
naturally acts transitively, whereas $\mathcal C_{\lambda}^{\,\rm{orb}}$
is parametrized by a closed interval on which $U(1)$ acts trivially.

In the following three subsections we determine the normal bundles, for
each type of  real touching conics. These subsections are organized as follows:
 first we
explicitly calculate the intersection of the irreducible components and some curves. 
Consequently we get a function of
$\lambda$ (which will be written $h_i$). 
Second we show that the normal bundle in problem degenerates into 
$O\oplus O(2)$ precisely when
$\lambda$ is a critical point of this function. Finally we
 determine the critical points.

 The consequences of the results in these three subsections 
will be postponed until \S
\ref{ss-con}.

\subsection{The case of generic type}\label{ss-generic}

Suppose $\lambda\in I_2\cup I_4^-\cup I_4^+$ and take
$C_{\theta}\in\mathcal C_{\lambda}^{\,\rm{gen}}$.
First we  calculate the intersection of $C_{\theta}$ and
$l_{\infty}$, where $l_{\infty}$ is the real line defined by $y_0=y_1=0$
as before. Let
$x_2=y_2/y_3$ be a non-homogeneous coordinate on $l_{\infty}$ (around
$P_{\infty}=(0:0:0:1)$).

\begin{lemma}\label{lemma-int-a} The set $\{C_{\theta}\cap l_{\infty}\set
C_{\theta}\in \mathcal C_{\lambda}^{\rm{gen}} \}$  consists of  disjoint two
circles about
$P_{\infty}$ in $l_{\infty}$, whose radiuses 
(with respect to the coordinate $x_2$ above)
 are
 given by
$$
h_0(\lambda):=\frac{Q+\sqrt{Q^2-f}}{\sqrt{f}}\hspace{3mm}{\rm{and}}\hspace{3mm}
h_0(\lambda)^{-1}=\frac{Q-\sqrt{Q^2-f}}{\sqrt{f}}$$  respectively, where
we put
$Q=Q(\lambda,1)$ and $f=f(\lambda)$ as before. 
\end{lemma}

Note that  we have $Q^2-f>0$ and $Q>\sqrt{f}$ by Proposition \ref{prop-necessa},
and therefore $h_0>1>h_0^{-1}>0$ holds.
Moreover, $h_0$ and $h_0^{-1}$ are differentiable on
 $I_2\cup I_4^-\cup I_4^+$.

\vspace{2mm}
\noindent Proof. On $H_{\lambda}=\{(y_0:y_1:y_2)\}$, $l_{\infty}$ is defined by 
$y_1=0$. Therefore by (\ref{eqn-a-0}) we readily have 
\begin{equation}\label{eqn-intersec100}C_{\theta}\cap l_{\infty}=
\left\{x_2=\frac{-Q\pm\sqrt{Q^2-f}}{\sqrt{f}}\cdot
e^{-i\theta}\right\}.\end{equation} This directly implies the claim of
the lemma. 
\proofend

\vspace{2mm} By Proposition \ref{prop-inv-a}, $\Phi^{-1}(C_{\theta})$
consists of two irreducible components, both of which are real rational
curves. We denote these components $L^+_{\theta}$ and $L^-_{\theta}$,
although there is no canonical way of making a distinction of these two.
Again by Proposition \ref{prop-inv-a}, 
$L^+_{\theta}$ and $L^-_{\theta}$ respectively form 
disjoint families 
$$\mathcal L_{\lambda}^+=\{L^+_{\theta}\set\theta\in\mathbf
R\}\hspace{3mm} {\rm{and}}\hspace{3mm}\mathcal
L_{\lambda}^-=\{L_{\theta}^-\set\theta\in\mathbf R\}$$ of (real and
smooth)  rational curves on
$S_{\lambda}=\Phi^{-1}(H_{\lambda})$. These are real members 
of real pencils on $S_{\lambda}$ and each member has no real point by 
Proposition \ref{prop-a} (iii). 
 Because $U(1)$  acts also on the parameter spaces ($=S^1$) of
$\mathcal L^+_{\lambda}$ and
$\mathcal L^-_{\lambda}$, the normal bundles of $L^+_{\theta}$ and
$L^-_{\theta}$ inside
$Z$ are  independent of the choice of $\theta$. 
The following proposition plays a
key role in determining the normal bundle:

\begin{prop}\label{prop-nbc} For any $L\in \mathcal L^+_{\lambda}
\cup \mathcal L^-_{\lambda}$, the normal bundle of $L$ in
$Z$ is isomorphic to either $ O(1)^{\oplus 2}$ or $ O\oplus O(2)$. Further, the
latter holds iff
$\lambda$ is a critical point of
$h_0(\lambda)$ defined in Lemma \ref{lemma-int-a}.
\end{prop}

In particular, members of $\mathcal L^+_{\lambda}$ and 
$ \mathcal L^-_{\lambda}$ have the same normal bundle in $Z$.

\vspace{2mm}
\noindent Proof. By Proposition \ref{prop-inv-pl}, $L$ is contained in the
smooth surface
$S_{\lambda}=\Phi^{-1}(H_{\lambda})$ and therefore we have an exact
sequence
$0\ra N_{L/S_{\lambda}}\ra N_{L/Z}\ra N_{S_{\lambda}/Z}|_L\ra 0$. 
By  Proposition \ref{prop-inv-a} (ii), we have $N_{L/S_{\lambda}}\simeq 
O_L$.
On the other hand, again by Proposition \ref{prop-inv-pl},
$S_{\lambda}$ is a smooth member of $|(-1/2)K_Z|$.
Therefore by adjunction formula we have
$K_S\simeq K_Z|_{S_{\lambda}}\otimes N_{S_{\lambda}/Z}\simeq
K_Z |_{S_{\lambda}}\otimes (-1/2)K_Z|_{S_{\lambda}}$ and hence 
$N_{S_{\lambda}/Z}\simeq (-1/2)K_Z|_{S_{\lambda}}\simeq -K_{S_{\lambda}}$.
Hence we get $N_{S_{\lambda}/Z}|_L\simeq-K_{S_{\lambda}}|_L\simeq
-K_L\otimes N_{L/S}\simeq O_L(2)$.
Therefore by the short exact sequence above, $N_{\lambda}:=N_{L/Z}$ is
isomorphic to either
$ O\oplus O(2)$ or $O(1)^{\oplus 2}$. Thus we get the first claim of
the proposition.
 
In order to show the second claim,
 we first explain the natural real structure on $\Gamma(N_{\lambda})$,
the space of sections of $N_{\lambda}$.
Since $L$ is real,  $\sigma$ naturally acts on $\Gamma(N_{\lambda})$ as the
complex conjugation. For $s\in\Gamma(N_{\lambda})$ we denote by ${\rm Re} s$ 
and ${\rm Im} s$ the real part and the imaginary part of $s$ respectively.
Namely, ${\rm Re} s=(s+\ol{\sigma(s)})/2$ and 
${\rm Im} s=(s-\ol{\sigma(s)})/2$.
Next recall that any one-parameter family of holomorphic  deformation of
$L$ in $Z$ naturally gives rise to a holomorphic section of
$N_{\lambda}$.  
 We have the following two one-parameter families of deformations of
$L$ in
$Z$: one is obtained by moving $L$ by $\mathbf C^*$-action,
where the $\mathbf C^*$-action is
the complexification of the $U(1)$-action. The other is
obtained by moving the parameter $\lambda$ in $\mathbf C$, while fixing 
$\theta$. Let $s\in\Gamma(N_{\lambda})$ and $t\in\Gamma(N_{\lambda})$ be
the holomorphic sections associated to the former and the latter  family
respectively. These are clearly linearly independent sections. Because the
$\mathbf C^*$-action preserves
$S_{\lambda}$, it follows from Proposition \ref{prop-inv-a} (ii) and (iii) that
each of the curves of the former family are disjoint. This implies that $s$ is
nowhere vanishing.  Next we consider the latter family. First, noting
$H_{\lambda}\cap H_{\lambda'}=l_{\infty}$ for $\lambda\neq\lambda'$,
 $t$ can be zero only on $\Phi^{-1}(l_{\infty})$. Suppose
$h_0'(\lambda)=0$, where the derivative is with respect to $\lambda\in\mathbf R$.
 Then we have
$(h_0^{-1})'(\lambda)=0$. Then, since $h$ is a holomorphic function of $\lambda$,
it can be easily derived by using the Cauchy-Riemann equation that
$\partial h/\partial \lambda=0$. In the same way, we have
$\partial h^{-1}/\partial \lambda=0$.
 These imply that $t$ vanishes on
$\Phi^{-1}(l_{\infty})\cap L=:\{z_{\lambda},\ol{z}_{\lambda}\}$.
  Therefore by Lemma \ref{lemma-nb}, we get
$N_{\lambda}\simeq O\oplus O(2)$. 

Next suppose $h_0'(\lambda)\neq 0$, so
that  
$(h_0^{-1})'(\lambda)\neq 0$. We claim that Re$(as+t)$ cannot vanish
at $z_{\lambda}$ and $\ol{z}_{\lambda}$ simultaneously,
for any $a\in\mathbf C$. 
Because $C_{\theta}$ intersects
$l_{\infty}$ transversally, $t$ also becomes a nowhere vanishing section
under our assumption. Hence the claim is true for $a=0$.
Putting $a=a_1+ia_2$, $a_1,a_2\in\mathbf R$,
we easily get
\begin{equation}\label{eqn-reim}
{\rm{Re}} (as+t)=a_1{\rm{Re}}s+({\rm{Re}} t-a_2{\rm{Im}} s).
\end{equation}
Since $s$ comes from the
$\mathbf C^*$-action, and since its real part corresponds to the 
$U(1)$-action,
$({\rm{Re}}s)(z_{\lambda})$ is  represented by the
tangent vector of the $U(1)$-orbit going through $z_{\lambda}$.
On the other hand, (\ref{eqn-intersec100}) implies that  
$({\rm{Re}}t)(z_{\lambda})$ is represented by a tangent vector
which is parallel to $({\rm{Im}}s)(z_{\lambda})$.
Hence by (\ref{eqn-reim}), we can deduce that 
${\rm{Re}} (as+t)(z_{\lambda})=0$ implies $a_1=0$ and 
\begin{equation}\label{eqn-reim2}
({\rm{Re}}t)(z_{\lambda})=a_2({\rm{Im}} s)(z_{\lambda}).
\end{equation}
Similarly,  ${\rm{Re}} (as+t)(\ol{z}_{\lambda})=0$ implies  $a_1=0$ and 
\begin{equation}\label{eqn-reim3}
({\rm{Re}}t)(\ol{z}_{\lambda})=a_2({\rm{Im}} s)(\ol{z}_{\lambda}).
\end{equation}
Suppose $a_2>0$. Then since $\{{\rm{Re}}s(z_{\lambda}),
{\rm{Im}}s(z_{\lambda})\}$ is an oriented basis of 
$T_{z_{\lambda}}(\Phi^{-1}(l_{\infty}))$
from the beginning,
(\ref{eqn-reim2}) implies that $\{{\rm{Re}}s(z_{\lambda}),
{\rm{Re}}t(z_{\lambda})\}$ is an oriented basis of
$T_{z_{\lambda}}(\Phi^{-1}(l_{\infty}))$. Further, we have
$$({\rm{Re}}s)(\ol{z}_{\lambda})=\sigma_*(({\rm{Re}}s)(z_{\lambda}))
\hspace{3mm}{\mbox{and}}\hspace{3mm}
({\rm{Re}}t)(\ol{z}_{\lambda})=\sigma_*(({\rm{Re}}t)(z_{\lambda})).$$
Hence we get by (\ref{eqn-reim3})
 $$({\rm{Im}} s)(\ol{z}_{\lambda})
=\frac{1}{a_2}({\rm{Re}}t)(\ol{z}_{\lambda})=
\frac{1}{a_2}\sigma_*\left(({\rm{Re}}t)(z_{\lambda})\right).$$
So we have $$\{({\rm{Re}}s)(\ol{z}_{\lambda}),({\rm{Im}} s)(\ol{z}_{\lambda})\}=
\{\sigma_*(({\rm{Re}}s)(z_{\lambda})),
 \sigma_*(({\rm{Re}}t)(z_{\lambda}))/a_2\}.$$
But since $\sigma$ is anti-holomorphic, $\sigma$ is orientation reversing.
Further, as is already seen, $\{{\rm{Re}}s(z_{\lambda}),
{\rm{Re}}t(z_{\lambda})\}$ is an oriented basis of  $T_{z_{\lambda}}
(\Phi^{-1}(l_{\infty}))$
 (if $a_2>0$). This implies that $\{\sigma_*(({\rm{Re}}s)(z_{\lambda})),
 \sigma_*(({\rm{Re}}t)(z_{\lambda}))/a_2\}$ is an anti-oriented basis
of  $T_{\ol{z}_{\lambda}}(\Phi^{-1}(l_{\infty}))$.  This contradicts to the fact that
$\{({\rm{Re}}s)(\ol{z}_{\lambda}),$ $({\rm{Im}} s)(\ol{z}_{\lambda})\}$ is an
oriented basis of $T_{\ol{z}_{\lambda}}(\Phi^{-1}(l_{\infty}))$.
Therefore, Re$(as+t)$ cannot vanish
at $z_{\lambda}$ and $\ol{z}_{\lambda}$ simultaneously, provided $a_2>0$.
Parallel arguments show the same claim holds for the case $a_2<0$.
Thus we have shown that  Re$(as+t)$ cannot vanish
at $z_{\lambda}$ and $\ol{z}_{\lambda}$ at the same time, as claimed.
On the other hand, it is obvious that 
$as+t$ does not vanish except $\{z_{\lambda},\ol{z}_{\lambda}\}$.
Therefore, the zero locus of $as+t$ consists of at most one point
for any $a\in\mathbf C$.
Since $s$ is a nowhere vanishing section, Lemma \ref{lemma-nb}
implies $N_{\lambda}\simeq O(1)^{\oplus 2}$.
\proofend

\begin{lemma}\label{lemma-bhv-1} Let $h_0=h_0(\lambda)$ be the positive
valued function on $I_2\cup I_4$ defined in Lemma \ref{lemma-int-a},
which is differentiable on
$I_2\cup I_4\backslash\{\lambda_0\}$. Then 
$h_0$  has a unique critical point on $I_2$, and has no critical point on
$I_4\backslash\{\lambda_0\}$. 
\end{lemma}

\noindent Proof. 
We have $Q(-1)>0$, $Q(0)>0$ and $Q(b/a)>0$ by
Proposition \ref{prop-realpoint} (i), and  
\begin{equation}\label{eqn-g-1}
h_0=\frac{Q+\sqrt{Q^2-f}}{\sqrt{f}}=
\frac{Q}{\sqrt{f}}+\sqrt{\frac{Q^2}{f}-1}.
\end{equation} From these,  it follows that 
$\lim_{\lambda\downarrow -1}h_0(\lambda)=\lim_{\lambda\uparrow
0}h_0(\lambda)=+\infty$.
Therefor $h_0$ has at least one critical
point on $I_2$, since $h_0$ is differentiable on $I_2$. So to prove the lemma it
suffices to show that this is a unique critical point on $I_2\cup
I_4\backslash\{\lambda_0\}$.

We consider the real valued function
$\gamma:=Q^2/f$ defined on
$I:=I_1\cup I_2\cup I_3\cup I_4$, which is clearly differentiable
on $I$.
Then $h_0=\sqrt{\gamma}+\sqrt{\gamma-1}$ on $I_2\cup I_4$, and we have 
$$h_0'=\gamma'\cdot\left(\frac{1}{2\sqrt{\gamma}}+
\frac{1}{2\sqrt{\gamma-1}}\right),$$
provided
$\lambda\neq\lambda_0$. Therefore on $I_2\cup
I_4\backslash\{\lambda_0\}$,  $h_0'(\lambda)=0$ iff
$\gamma'(\lambda)=0$. 
It is readily seen  that
$\lim_{\lambda\da-\infty}\gamma(\lambda)=
\lim_{\lambda\ua-1}\gamma(\lambda)=-\infty$,
$\lim_{\lambda\da-1}\gamma(\lambda)=\lim_{\lambda\ua0}\gamma(\lambda)=\infty$,
$\lim_{\lambda\da 0}\gamma(\lambda)=\lim_{\lambda\ua
(b/a)}\gamma(\lambda)=-\infty$, and
$\lim_{\lambda\da
(b/a)}\gamma(\lambda)=\lim_{\lambda\ua\infty}\gamma(\lambda)=\infty$.  
Therefore
$\gamma$ has at least one critical point on each $I_j$, $1\leq j\leq 4$.
We also have  
\begin{equation}\label{eqn-diff}
 \gamma'=Q(2Q'f-Qf')/f^2.
\end{equation}
Suppose that  critical points of  $\gamma$ on $I_2$ are not unique. Then
$\gamma$ has at least three critical points on $I_2$.
This implies that $\gamma$ has at least four critical points on $I_2\cup I_4$. 
Because $Q>0$  on $I_2\cup I_4$ (Proposition
\ref{prop-realpoint} (i)), 
these critical points must correspond to zeros of
$2Q'f-Qf'$ whose degree is just four.
By (\ref{eqn-diff})
this implies that the other  critical points of $\gamma$
on
$I_1$ and $I_3$ must correspond to zeros of $Q$.
But this cannot happen since $Q>0$  on
$I_2\cup I_4$ and since $Q$ is degree two.
Therefore, our assumption fails and it follows that critical points on $\gamma$ 
on $I_2$ is unique. Hence  critical points of $h_0$ on $I_2$ is also
unique. Exactly the same argument shows that
$\gamma$  has a unique critical point on $I_4$.
This critical point must be $\lambda_0$, since 
$\gamma$ attains the minimal value $(=1)$ there. This implies that
$g$ has no critical point on $I_4\backslash\{\lambda_0\}$. Thus we
obtain the claims of the lemma.
\proofend

\vspace{2mm} The following is the main result of this subsection: 

\begin{prop}\label{prop-nottl} (i) If $\lambda\in I_4$ and if
$\lambda\neq
\lambda_0$, we  have $N_{L/Z}\simeq O(1)^{\oplus 2}$ for any
$L\in\mathcal L_{\lambda}^+\cup\mathcal L_{\lambda}^-$. (ii) There is a
unique
$\lambda\in I_2$ such that
$N_{L/Z}\simeq O\oplus O(2)$ for any $L\in\mathcal
L_{\lambda}^+\cup\mathcal L_{\lambda}^-$. For any other $\lambda\in
I_2$,  we  have $N_{L/Z}\simeq O(1)^{\oplus 2}$ for arbitrary $L\in\mathcal
L_{\lambda}^+\cup\mathcal L_{\lambda}^-$.  (iii) If $\lambda\in I_2$, 
any  member of $\mathcal L_{\lambda}^+\cup\mathcal L_{\lambda}^-$ is not
a twistor line in $Z$ (even if
$Z$ is actually a twistor space).
\end{prop}

\noindent Proof. (i) and (ii) are direct consequences of Proposition
\ref{prop-nbc} and Lemma
\ref{lemma-bhv-1}.  To show (iii), let $\lambda'\in I_2$ be the unique
critical point of $g$. Then by (ii), any $L'\in\mathcal
L_{\lambda'}^+\cup\mathcal L_{\lambda'}^- $ is not a twistor line. 
We can see that for any
$\lambda\in I_2$ and for any $L\in\mathcal L_{\lambda}^+\cup\mathcal
L_{\lambda}^-$, $L$ can be deformed into some $L'\in\mathcal
L_{\lambda'}^+\cup\mathcal L_{\lambda'}^- $ preserving the real
structure. In fact, we have $\Phi(L)=C_{\theta}$ for some
$C_{\theta}\in\mathcal C_{\lambda}^{\,\rm{gen}}$.  Then since $I_2$ is an interval in
$\mathbf R$,
$C_{\theta}$ can be canonically deformed into some
$C_{\theta}'\in \mathcal C_{\lambda'}^{\,\rm{gen}}$. (The point is that we take
a constant
$\theta$ for any $\lambda\in I_2$.)  Correspondingly, we obtain
deformation of
$L$ into  $L'\in\mathcal L_{\lambda'}^+\cup\mathcal L_{\lambda'}^- $
such that
$\Phi(L')=C'_{\theta}$. Thus we get an explicit real one-dimensional
 family of  rational curves in
$Z$ containing $L$ and $L'$ as its members,
 as claimed. Any member of this family is
 real by  Proposition \ref{prop-inv-a} (iv).  Since any
deformation of twistor line preserving the real structure is still a
twistor line, it follows that
$L$ is not a twistor line.
\proofend
 
\vspace{2mm} We note the proof of (iii) does not work for $I_4$, since as
$\lambda$ goes to
$\lambda_0$, the curve
$C_{\theta}$ (defined by (\ref{eqn-a-0})) degenerates into a double line.
This is an important point for obtaining a natural compactification of
the space of real touching conics of $B$.

\begin{cor}\label{cor-image1} If $\lambda\in I_2$, only the members of
$\mathcal C_{\lambda}^{\,\rm{orb}}$ can be the image of twistor lines.
Namely  over $I_2$, members of $\mathcal C_{\lambda}^{\,\rm{gen}}$ cannot be the
images of twistor lines.
\end{cor}

\noindent Proof.
By Proposition \ref{prop-type} (i), the image of a twistor line contained
in
$H_{\lambda}$ is either a touching conic of generic type or that of orbit
type for $\lambda\in I_2\cup I_4$. But by Proposition \ref{prop-nottl}
(iii) the former cannot be the image of a twistor line if $\lambda\in
I_2$.
\proofend

\begin{cor}\label{cor-chatl} A twistor line of a self-dual 4-manifold (i.e.\! a
fiber of the twistor fibration) is not in general characterized by  the property
that it is a real smooth rational curve without fixed point whose normal bundle is
isomorphic to $O(1)^{\oplus 2}$. 
More concretely, the twistor space of any non-LeBrun self-dual metric on
$3\mathbf{CP}^2$ of positive scalar curvature with a non-trivial 
Killing field always possesses such a real rational curve.
\end{cor}

\noindent Proof. 
Let $Z$ be a twistor space as in the corollary.
Then $Z$ has a
structure as in  Proposition \ref{prop-def-B}, where $Q,a$ and $b$
satisfy the conditions in Proposition \ref{prop-necessa}. By (ii) and (iii)
of Proposition
\ref{prop-nottl}, $Z$ always has a real smooth rational
curve $L$ satisfying
$N_{L/Z}\simeq O(1)^{\oplus 2}$, but which is not a twistor line.
This $L$ has no real point by Proposition \ref{prop-a} (iii) and
the reality of $\Phi$.
Existence of $Z$ is proved in \cite{Hon02} and \cite{Hon-II}.
 \proofend

\vspace{2mm} We remark that C.\,Simpson \cite{S97}  asked a similar
question about a characterization of twistor lines for
 twistor spaces of  hyperK\"ahler manifolds.  

Next we give another geometric proof for the fact that $L$ cannot be a
twistor line for
$\lambda\in I_2$ (although we will not need this result in the sequel).

\begin{prop} \label{prop-break} If
$\lambda\in I_2$  is not a critical point of $h_0$, there exists a unique
$\mu\in I_2$ with
$\lambda\neq\mu$ satisfying the following: for any $L\in\mathcal
L_{\lambda}^+$ (resp.\,$L\in\mathcal L^-_{\lambda}$) there exists
$L'\in\mathcal L_{\mu}^+$ (resp.\,$L'\in\mathcal L^-_{\mu}$) such that $L\cap
L'\neq\phi$. 
\end{prop}

\noindent Proof. Let $\lambda'\in I_2$ be the unique critical point of
$h_0$ as before. By our proof of Lemma \ref{lemma-bhv-1} we have
$\lim_{\lambda\downarrow -1}h_0(\lambda)=\lim_{\lambda\uparrow
0}h_0(\lambda)=+\infty$ and
$g$ is strictly decreasing on $(-1,\lambda')$ and strictly
increasing on
$(\lambda',0)$. Suppose $\lambda<\lambda'$. Let $l$ and $\ol{l}$ be the 
conjugate pair of rational curves which are mapped biholomorphically onto
$l_{\infty}$. (See the proof of Proposition \ref{prop-image}.)  Then by
Lemma
\ref{lemma-int-a}, $L\cap l$ is a point  which is  either
$x_2=h_0(\lambda) e^{i\theta}$ or
$x_2=h_0(\lambda)^{-1}e^{i\theta}$ for some 
$\theta\in\mathbf R$,  where we identify $l$ and
$l_{\infty}$ via $\Phi$ and use  $x_2=y_2/y_3$  as an affine coordinate on
$l_{\infty}$ as before. If $x_2=h_0(\lambda)^{-1}e^{i\theta}$, 
$\ol{l}\cap L$ is a point having $x_2=h_0(\lambda)e^{i\theta}$. Thus by
a possible exchange of $l$ and
$\ol{l}$, we may suppose that $l\cap L$ is a point satisfying
$x_2=h_0(\lambda)e^{i\theta}$. Then by the behavior of $h_0$ mentioned above,
 there exists a unique
$\mu>\lambda'$, $\mu\in I_2$  such that $h_0(\lambda)=h_0(\mu)$. On the
other hand, by our choice of $L$ we have $\Phi(L)=C_{\theta}$ for some
$C_{\theta}\in\mathcal C_{\lambda}^{\,\rm{orb}}$. Then take $L'\in\mathcal L_{\mu}^+$
such that
$\Phi(L')=C_{\theta}\in \mathcal C_{\mu}^{\,\rm{orb}}$. (Although we use the same symbol
$C_{\theta}$, they represent different conics since $\lambda\neq \mu$.
The point is that we take the same $\theta$ for different $\lambda$'s.)
Then
$L\cap L'\cap l$ is  a point satisfying
$x_2=h_0(\lambda)e^{i\theta}$. 
(We also have $L\cap L'\cap \ol{l}$ is a point satisfying
 $x_2=h_0(\lambda)^{-1}e^{i\theta}$.)
Thus we have proved the claim for
$\lambda<\lambda'$. Of course, the case $\lambda>\lambda'$
and the case $L\in \mathcal L^{-}_{\lambda}$ are similar.
\proofend
 
\vspace{2mm} The proposition shows that when $\lambda\in I_2$ passes
through the critical point ($=\lambda'$) of
$h_0$,
 the local  twistor fibration arising from $L\in\mathcal L_{\lambda}^+\cup
\mathcal L^-_{\lambda}$, $\lambda\neq \lambda'$ breaks down. Note also
that Proposition
\ref{prop-break} also holds true for $I_4$ without any change of the
proof, and it implies that if members of $\mathcal L_{\lambda}^+$ (resp.
\!$\mathcal L_{\lambda}^-$) are twistor lines for $\lambda\in I_4^-$, 
 members of $\mathcal L_{\lambda}^-$ (resp. $\mathcal L_{\lambda}^+$)
must be  twistor lines for $\lambda \in I_4^+$.

\subsection{The case of special type}\label{ss-special} In this subsection we
calculate the normal bundle of $L^+$ and $L^-$ in $Z$, where
$L^+$ and $L^-$ are  curves which are mapped biholomorphically onto a real
touching conic of special type. Compared to the case of generic type, the
problem becomes  harder and the result becomes  more complicated, since
touching conics of special type go through the singular point
$P_{\infty}$ and
$\ol{P}_{\infty}$ of
$B$, so that  the situation, and hence  the result also,  depend on how
we resolve the corresponding singularities of $Z_0$.

First we recall the situation and fix notations. Let $\Phi_0:Z_0\ra
\mathbf{CP}^3$ be the double covering branched along $B$. Put
$p_{\infty}:=\Phi_0^{-1}(P_{\infty})$.
 In a neighborhood of $P_{\infty}=(0:0:0:1)$, we use 
$(x_0,x_1,x_2)$ as an affine coordinate by setting $x_i=y_i/y_3$. Then around
$P_{\infty}=(0,0,0)$,
$B$ is given by the equation
$(x_2+Q(x_0,x_1))^2-x_0x_1(x_0+x_1)(ax_0-bx_1)=0$.  Let $z$ be a fiber
coordinate of $O(2)\ra\mathbf{CP}^3$. Then $Z_0$ is given by the
equation 
\begin{equation}\label{eqn-Z_0-2}
z^2+\left(x_2+Q(x_0,x_1)\right)^2-x_0x_1(x_0+x_1)(ax_0-bx_1)=0.
\end{equation} This can be also written as
$\{z+i(x_2+Q(x_0,x_1))\}\{z-i(x_2+Q(x_0,x_1))\}=x_0x_1(x_0+x_1)(ax_0-bx_1)$.
Setting $\xi=z+i(x_2+Q(x_0,x_1))$ and $\eta=z-i(x_2+Q(x_0,x_1))$, we get
\begin{equation} Z_0:\hspace{3mm}\xi\eta=x_0x_1(x_0+x_1)(ax_0-bx_1).
\end{equation} Thus $p_{\infty}=\{(x_0,x_1,\xi,\eta)=(0,0,0,0)\}$ is a compound
$A_3$-singularity. Small resolutions of $p_{\infty}$ are explicitly
given as follows: first we choose ordered three linear forms
$\{\ell_1,\ell_2,\ell_3\}\subset\{x_0,x_1,x_0+x_1,ax_0-bx_1\}$. Next
blow-up
$Z_0$ along
$\{\xi=\ell_1=0\}$. Then by setting $\xi=u\ell_1$ we get 
\begin{equation} Z'_0: u\eta=x_0x_1(x_0+x_1)(ax_0-bx_1)/\ell_1.
\end{equation} The exceptional curve of $Z'_0\ra Z_0$ is given by
$\Gamma_1':=\{(u,\eta,x_0,x_1)\set \eta=x_0=x_1=0\}$. We can use $u$
as an affine coordinate on $\Gamma_1'$.
 Next blowing up $Z'_0$  along $\{u=\ell_2=0\}$  and setting
$u=v\ell_2$, we get
\begin{equation} Z''_0: v\eta=x_0x_1(x_0+x_1)(ax_0-bx_1)/\ell_1\ell_2.
\end{equation} The exceptional curve of $Z''_0\ra Z'_0$ is given by
$\Gamma_2'':=\{(v,\eta,x_0,x_1)\set \eta=x_0=x_1=0\}$, on which we can
use $v$ as an affine coordinate.
 Finally by blowing up along $\{v=\ell_3=0\}$ and setting $v=w\ell_3$,
we get 
\begin{equation} Z: w\eta=x_0x_1(x_0+x_1)(ax_0-bx_1)/\ell_1\ell_2\ell_3,
\end{equation} which is clearly smooth in a neighborhood of the origin. The
exceptional curve of $Z\ra Z''_0$ is given by
$\Gamma_3:=\{(w,\eta,x_0,x_1)\set
 \eta=x_0=x_1=0\}$, on which we can use $w$ as an affine coordinate. 


Once a resolution of $p_{\infty}$ is given, it naturally determines that of 
$\ol{p}_{\infty}$
by reality.
Let $\mu$ be the resolution of $p_{\infty}$ and $\ol{p}_{\infty}$
(for some choice of $\ell_1$, $\ell_2$ and $\ell_3$).
Thus a choice of $\ell_1,\ell_2$ and
$\ell_3$ determines a small resolution of $p_{\infty}$
and $\ol{p}_{\infty}$ and there are
$4!=24$ ways of resolutions in all.
Let $\Gamma=\mu^{-1}(p_{\infty})=\Gamma_1+\Gamma_2+\Gamma_3$ and
$\ol{\Gamma}=\mu^{-1}(\ol{p}_{\infty})=
\ol{\Gamma}_1+\ol{\Gamma}_2+\ol{\Gamma}_3$ be the
exceptional curves  of $\mu$, where $\Gamma_i$ and $\ol{\Gamma}_i$
are the exceptional
curves arising from the
$i$-th blowing-up above. Then we have
$\Gamma_1\cap\Gamma_2\neq\phi,\Gamma_2\cap\Gamma_3\neq\phi$ and
$\Gamma_1\cap\Gamma_3=\phi$. 

In order to calculate the intersection $L^+$ and $L^-$ with $\Gamma$,
we  need a one-parameter presentation of $C_{\theta}$, in a neighborhood of
$p_{\infty}$:

\begin{lemma}\label{lemma-one-rep} Let $C_{\theta}\subset H_{\lambda}$
be a real touching conic of special type whose equation is given by
(\ref{eqn-b-0}), and $(x_0,x_1,x_2)$ the affine coordinate around
$P_{\infty}$ as above. Then in  a neighborhood of
$P_{\infty}$,
$C_{\theta}$ has a one-parameter presentation of the following form:
\begin{equation}\label{eqn-solve}
\left\{
\begin{array}{l}x_0=\lambda x_1\\ \displaystyle
x_2=-Be^{-i\theta}x_1-\frac{\sqrt{Q^2-f}+Q}{2}x_1^2+
\frac{\sqrt{Q^2-f}+Q}{2}Be^{i\theta}x_1^3+O(x_1^4),
\end{array}\right.
\end{equation}
where we put
$$B:=B(\lambda)=\left(\frac{\sqrt{Q^2-f}-Q}{2}\right)^{\frac{1}{2}}.$$
\end{lemma}

Note again that $f<0$ guarantees $\sqrt{Q^2-f}-Q>0$.

\vspace{2mm}
\noindent Proof. By  solving (\ref{eqn-b-0}) with respect to $x_2$, we
get
\begin{equation}\label{eqn-solv} x_2=-g(x_1)\cdot x_1,\hspace{3mm}
g(x_1):=
\frac{Be^{-i\theta}+\sqrt{Q^2-f}\,\,x_1}{1+Be^{i\theta}x_1}.
\end{equation} Calculating the Maclaurin expansion of $g(x_1)$, we get
(\ref{eqn-solve}). This is a routine work and we omit the
detail.\proofend

\begin{lemma}\label{lemma-tlr} In a neighborhood of $p_{\infty}$, each of the two
irreducible components of
$\Phi^{-1}_0(C_{\theta})$ has a one-parameter presentation with
respect to
$x_1$ in the following forms respectively:
\begin{equation}\label{eqn-L1}
\xi=-2iBe^{-i\theta}x_1+O(x_1^2),\hspace{3mm}
\eta=\frac{ie^{i\theta}f}{2B}x_1^3+O(x_1^4),\hspace{3mm} x_0=\lambda x_1,
\end{equation} and
\begin{equation}\label{eqn-L2}
\xi=-\frac{ie^{i\theta}f}{2B}x_1^3+O(x_1^4),\hspace{3mm} 
\eta=2iBx_1e^{-i\theta}+O(x_1^2),\hspace{3mm} x_0=\lambda x_1.
\end{equation}

\end{lemma}

\noindent Proof. First by substituting $x_0=\lambda x_1$ into
(\ref{eqn-Z_0-2}), we get
$$ z^2=(f-Q^2)x_1^4-2Qx_1^2x_2-x_2^2.
$$ Substituting (\ref{eqn-solv}) into this, we get 
$$z^2=\left\{(f-Q^2)x_1^2+2Qg(x_1)x_1-g(x_1)^2\right\}\,x_1^2.$$ Hence we have
$$z=\pm k(x_1)\,x_1,\hspace{3mm}
k(x_1)=\left\{(f-Q^2)x_1^2+2Qg(x_1)x_1-g(x_1)^2
\right\}^{\frac{1}{2}}.$$ From this
we deduce 
$$
\xi=z+i(x_2+Qx_1^2)=\left(\pm k(x_1)-ig(x_1)\right)\,x_1+iQ x_1^2
$$  and 
$$\eta=z-i(x_2+Qx_1^2)=\left(\pm k(x_1)+ig(x_1)\right)\,x_1-iQ x_1^2.$$ Then we get
the  desired equations by calculating the Maclaurin expansions of 
$\pm k(x_1)-ig(x_1)$ and $\pm k(x_1)+ig(x_1)$. These are also  routine
works and we omit the detail.
\proofend

\vspace{2mm}

\begin{lemma}\label{lemma-interse} Let $L^+_{\theta}$ and
$L^-_{\theta}$  be the curves in $Z$ which are the proper transforms of the
curves (\ref{eqn-L1}) and (\ref{eqn-L2}) respectively.  Then we have: (i)
$L^+_{\theta}\cap\Gamma_1$ is a point satisfying 
\begin{equation}\label{eqn-inters}
u=-2iBe^{-i\theta}\cdot\frac{x_1}{\ell_1}
\end{equation}
and $L^+_{\theta}\cap
\Gamma_2$ and $L^+_{\theta}\cap\Gamma_3$ are empty, (ii)
$L^-_{\theta}\cap\Gamma_1$ and $L^-_{\theta}\cap\Gamma_2$ are empty and
$L^-_{\theta}\cap\Gamma_3$ is a point satisfying
\begin{equation}\label{eqn-inter2}
w=-\frac{ie^{i\theta}f}{2B}\cdot
\frac{x_1^3}{\ell_1\ell_2\ell_3}
\end{equation} 
where we use $u$ and
$w$ as local coordinates on $\Gamma_1$ and $\Gamma_3$ respectively as
explained before, and $B=B(\lambda)$ is as in Lemma
\ref{lemma-one-rep}.  
\end{lemma}

Here note that
$x_1/\ell_1$ and $x_1^3/\ell_1\ell_2\ell_3$ do not depend on $x_1$, and
depend  on $\lambda$ only. Further, we have $B>0$ since $f<0$.

\vspace{2mm}
\noindent Proof. By substituting $\xi=u\ell_1$ into (\ref{eqn-L1}), we
get the inverse image of (\ref{eqn-L1}) to be
$$ u\ell_1=-2iBe^{-i\theta}x_1+O(x_1^2),\hspace{3mm}
\eta=\frac{ie^{i\theta}f}{2B}x_1^3+O(x_1^4),\hspace{3mm} x_0=\lambda x_1.
$$ Eliminating  $\Gamma_1'=\{
\eta=x_0=x_1=0\}$, we get the equation of the
proper transform in $Z'_0$ to be
$$ u=-2iBe^{-i\theta}\frac{x_1}{\ell_1}+O(x_1),\hspace{3mm}
\eta=\frac{ie^{i\theta}f}{2B}x_1^3+O(x_1^4),\hspace{3mm} x_0=\lambda x_1.
$$ By setting $x_1=0$, we get $u=-2iBe^{-i\theta}\cdot x_1/\ell_1$. As
remarked above, $B$ is non-zero. Therefore the remaining blow-ups 
$Z''_0\ra Z_0'$ and $Z\ra
Z''_0$ do not have effect on the intersection. Hence we get
(\ref{eqn-inters}). Similar calculations show (\ref{eqn-inter2}). (Note
that we have
$\xi=w\,\ell_1\ell_2\ell_3$.)
 \proofend

\vspace{2mm} As in the previous subsection we put $\mathcal
L_{\lambda}^+=\{L^+_{\theta}\set\theta\in\mathbf R\}$ and $\mathcal
L_{\lambda}^-=\{L^-_{\theta}\set\theta\in\mathbf R\}$. (Note that this
time we explicitly specified $L^+_{\theta}$ and $L^-_{\theta}$
respectively in Lemma \ref{lemma-tlr}.)
$U(1)$ again acts transitively  on the parameter spaces of these
families.  By Lemma
\ref{lemma-interse},
$L^+_{\theta}\cap \Gamma_1$ is a point, and 
$\{L^+_{\theta}\cap \Gamma_1\set\theta\in\mathbf R\}$ is a circle in
$\Gamma_1$ whose radius is
$$h_1(\lambda):=2B\cdot |x_1/\ell_1|.$$ Similarly,
$\{L^-_{\theta}\cap
\Gamma_3\set\theta\in\mathbf R\}$ is a circle in $\Gamma_3$ whose radius
is
$$h_3(\lambda):= (-f/2B)\cdot |x_1^3/\ell_1\ell_2\ell_3|.$$ 
($h_2$ will appear in the next subsection.) Then  we
have the following proposition, which implies that the normal bundles of
$L^+_{\theta}$ and
$L^-_{\theta}$  in $Z$ are  determined by the behavior of
$h_1$ and $h_3$ respectively.

\begin{prop}\label{prop-nbc2} Assume $\lambda\in I_1\cup I_3$ and take
any  $L=L^+_{\theta}\in\mathcal L_{\lambda}^+$. Then either
$N_{L/Z}\simeq O(1)^{\oplus 2}$ or
$N_{L/Z}\simeq O\oplus O(2)$ holds. Further, the latter holds iff
$\lambda$ is a critical point of
$h_1$ above. The same claims  also hold  for $\mathcal L^-_{\lambda}$ if
we replace $h_1$ by $h_3$. 
\end{prop}

\noindent Proof. The first claim can be proved in the same way as in 
Proposition \ref{prop-nbc}.
The other claims can also be  proved in the same manner as
in Proposition \ref{prop-nbc}:
take any $L\in \mathcal L_{\lambda}^+$. Then the two
one-parameter families of $L$ in $Z$  in the previous proof make senses
also in this case, so that  we again have two linearly independent
sections $s$ and $t$ of $N_{\lambda}$, $N_{\lambda}=N_{L/Z}$.
Then the previous proof works if we replace $h_0$ by $h_1$, 
$\Phi^{-1}(l_{\infty})$ by $\Gamma_1\cup \ol{\Gamma}_1$, and
(\ref{eqn-intersec100}) by (\ref{eqn-inters}).
For $L\in \mathcal L_{\lambda}^-$,
replace $h_0$ by $h_3$, 
$\Phi^{-1}(l_{\infty})$ by $\Gamma_3\cup \ol{\Gamma}_3$, and
(\ref{eqn-intersec100}) by (\ref{eqn-inter2}).
\proofend

\vspace{2mm} By definition, $h_1$ and $h_3$ depend on the choice of
$\ell_1,\ell_2$ and $\ell_3$. Therefore, Proposition \ref{prop-nbc2}
implies that 
 the normal bundles of $L^+$ and $L^-$  in $Z$ depend on how we resolve
$p_{\infty}$. More precisely, the normal bundle of $L^+$ depends on the
choice of $\ell_1$ only, whereas  the normal bundle  of $L^-$ depends on
that of
$\{\ell_1,\ell_2,\ell_3\}$.

Thus we need to know the critical points of $h_1$ and $h_3$ for every
choices of
$\ell_1,\ell_2$ and $\ell_3$.  At first glance there may seem to be too
many functions to be investigated, but it is easily seen that $h_3$ is a
reciprocal of $h_1$ (up to a constant) for some other choice of
$\ell_1,\ell_2$ and $\ell_3$. Therefore, what we need to know is the
behavior of $h_1$ for  the  four choices of $\ell_1$.
(Behavior of these functions near the endpoints of $I_1$ and $I_3$ below will
be needed  in \S \ref{ss-con}.)

\begin{lemma}\label{lemma-crtcl-sp} (i) If $\ell_1=x_1$,
$h_1$ has no critical point on
$I_1$, and has a unique critical point on $I_3$.  Further, we have
$\lim_{\lambda
\da-\infty}h_1(\lambda)=+\infty$ and $h_1(-1)=0$. (ii)  If $\ell_1=x_0$,
$h_1$ has a unique critical point on
$I_1$, and no critical point on $I_3$. Further, we have $\lim_{\lambda
\da 0}h_1(\lambda)=+\infty$ and $h_1(b/a)=0$. (iii) If
$\ell_1=x_0+x_1$,
$h_1$ has no critical point on
$I_1$, and has a unique critical point on $I_3$.  Further, we have
$\lim_{\lambda
\da-\infty}h_1(\lambda)=0$ and $\lim_{\lambda\uparrow
-1}h_1(\lambda)=+\infty$. (iv)  If $\ell_1=ax_0-bx_1$, $h_1$ has a
unique critical point on
$I_1$, and no critical point on $I_3$. Further, we have $h_1(0)=0$ and
$\lim_{\lambda\uparrow b/a}h_1(\lambda)=+\infty$.
\end{lemma}


\noindent Proof.  (i) If
$\ell_1=x_1$, we have $h_1^2=2(\sqrt{Q^2-f}-Q)$.
Since $h_1>0$ on $I_1\cup I_3$, the critical points of 
$h_1^2$ and $h_1$ coincide on $I_1\cup I_3$.
We think of $h^2_1$  as a real valued function defined on the whole
of $\mathbf R$, but which is not differentiable at $\lambda=\lambda_0$
in general. It is immediate to see that
$h_1^2(-1)=h_1^2(0)=h_1^2(b/a)=0$,
$\lim_{\lambda\da-\infty}h_1^2(\lambda)=+\infty$ and 
$\lim_{\lambda\ua\infty}h_1^2(\lambda)=-\infty$.
Then because $h_1^2$ is differentiable on $\lambda\neq\lambda_0$,
$h_1^2$ has a critical point on $I_2$ and $I_3$ respectively.
On the other hand, we have
$$
\left(\sqrt{Q^2-f}-Q\right)'=\frac{2QQ'-f'-2Q'\sqrt{Q^2-f}}{2\sqrt{Q^2-f}},
$$ and it follows that $(\sqrt{Q^2-f}-Q)'=0$ implies
\begin{equation}\label{eqn-diff00}
(2QQ'-f')^2=4Q'\
\!^2(Q^2-f).
\end{equation} It is readily seen that the degree of both hand sides
of (\ref{eqn-diff00}) are six, and that both have $\lambda_0$ as a double root.
have  $(\lambda-\lambda_0)^2$ as a factor.
Since we have already got two critical points of $h_1^2$
other than $\lambda=\lambda_0$, 
there are at most two solutions of (\ref{eqn-diff00})  remaining.
 
We set $g:=2(-\sqrt{Q^2-f}-Q)$ which is also defined on $\mathbf R$ and
possibly not differentiable at $\lambda=\lambda_0$.
Note that if we replace $h_1^2$ by $g$ on $\lambda\geq\lambda_0$,
then the resulting function is differentiable at $\lambda=\lambda_0$.
It is  easily verified that 
$g'=0$ also implies (\ref{eqn-diff00}) and it gives a solution
not coming from $(h_1^2)'=0$.
Further, we readily have 
$\lim_{\lambda\da-\infty}g(\lambda)=
\lim_{\lambda\ua\infty}g(\lambda)=-\infty$.

Suppose that $g$ has a critical point.
Then together with the above two critical points of $h_1^2$ on $I_2\cup
I_3$, we have three solutions of (\ref{eqn-diff00}) other than
$\lambda=\lambda_0$. If $h_1^2$ has  critical points on $I_1$, its
number is at least two. This implies that (\ref{eqn-diff00}) has five
solutions other than
$\lambda=\lambda_0$ and this is a contradiction.
Therefore $h_1^2$  has no critical points on $I_1$,
if $g$ has a critical point.
Similarly, if the number of the critical points on $I_3$ is not one,
then it must be at least three.
This is again
a contradiction. Thus if $g$ has a critical point,
$h_1^2$, and hence $h_1$ has no critical point on $I_1$ and a
unique critical point on $I_3$.

So suppose that $g$ has no critical point.
This happens exactly when $g$ attains the maximal
value at $\lambda=\lambda_0$. Then we have $\lim_{\lambda\ua\lambda_0}
g'(\lambda)>0$,
since otherwise $g$ has a critical point on $\lambda<\lambda_0$.
Because we have $\lim_{\lambda\ua\lambda_0}
g'(\lambda)=\lim_{\lambda\da\lambda_0}(h_1^2)'(\lambda)$,
we get $\lim_{\lambda\da\lambda_0}(h_1^2)'(\lambda)>0$.
Since $\lim_{\lambda\ua\infty}h_1^2(\lambda)=-\infty$,
it follows that $h_1^2$ has a critical point on $I_4$. Thus we get
three solutions of (\ref{eqn-diff00}) other than $\lambda_0$.
Then the same argument in the case that $g$ has a critical point as above,
we can deduce that $h_1$ has no critical point on $I_1$ and
a unique critical point on $I_3$.
Thus we get the claim of (i) concerning  critical points of $h_1$.
The remaining claims of (i) immediately follows from the definition of
$h_1$.

Claims of 
(ii), (iii) and (iv) about critical points can be obtained by applying a
projective transformation $\lambda\mapsto 1/\lambda$ for the case (ii), 
$\lambda\mapsto 1/(\lambda+1)$ for the case (iii), and 
$\lambda\mapsto 1/(a\lambda-b)$ for the case (iv) respectively.
The other claims are immediate to see.
\proofend

\vspace{2mm} As is already mentioned, the behavior of $h_3$ can be
easily seen from that of $h_1$ for some other choice of $\ell_1,\ell_2$
and $\ell_3$. The result is the following:

\begin{lemma}\label{lemma-crtcl-sp2} (i) If
$\{\ell_1,\ell_2,\ell_3\}=\{x_0,x_0+x_1,ax_0-bx_1\}$, $h_3$ has no
critical point on
$I_1$, and has a unique critical point on $I_3$.  Further, we have
$\lim_{\lambda
\da-\infty}h_3(\lambda)=0$ and $\lim_{\lambda\ua
-1}h_3(\lambda)=+\infty$. (ii)  If
$\{\ell_1,\ell_2,\ell_3\}=\{x_1,x_0+x_1,ax_0-bx_1\}$,
 $h_3$ has a unique critical point on
$I_1$, and no critical point on $I_3$. Further, we have $h_3(0)=0$ and
$\lim_{\lambda\ua b/a}h_3(\lambda)=+\infty$.   (iii) If
$\{\ell_1,\ell_2,\ell_3\}=\{x_0,x_1,ax_0-bx_1\}$, $h_3$ has no critical
point on
$I_1$, and has a unique critical point on $I_3$. 
 Further, we have $\lim_{\lambda
\da-\infty}h_3(\lambda)=+\infty$ and $h_3(-1)=0$. (iv)  If
$\{\ell_1,\ell_2,\ell_3\}=\{x_0,x_1,x_0+x_1\}$, $h_3$ has a unique
critical point on
$I_1$, and no critical point on $I_3$. Further, we have $\lim_{\lambda\da
0}h_3(\lambda)=+\infty$ and
$h_3(b/a)=0$.
\end{lemma}

\begin{cor}\label{cor-oneof} For any choice of $\ell_1,\ell_2$ and
$\ell_3$,  the following (i) and (ii) hold: (i) 
members of $\mathcal L^+_{\lambda}$ ($\lambda\in I_1$) and  $\mathcal
L^+_{\lambda}$ ($\lambda\in I_3$) cannot be twistor lines at the
same time,  (ii)
the same claim holds also for $\mathcal L^-_{\lambda}$.
\end{cor}

\noindent Proof. Suppose that $\lambda\in I_1\cup I_3$ is a critical
point of
$h_1$. Then by Proposition \ref{prop-nbc2}, any member of $\mathcal
L^+_{\lambda}$ is not a twistor line because its normal bundle in $Z$ is
$O\oplus O(2)$. Then just as in the proof of Proposition
\ref{prop-nottl}, any member of $\mathcal L^+_{\mu}$ cannot be a twistor
line provided that $\mu\in I_1\cup I_3$ and $\lambda$ belong to the same 
interval ($I_1$ or $I_3$). By Lemma \ref{lemma-crtcl-sp},
 $h_1$  necessarily  has  a critical point on just one of $I_1$ and
$I_3$. Hence (i) holds. The proof is the same for $\mathcal
L^-_{\lambda}$ if we use  Lemma \ref{lemma-crtcl-sp2} instead.
\proofend

\vspace{2mm} Thus together with Proposition \ref{prop-nbc2},
 we have obtained  new families of
real smooth rational curves which have
$O(1)^{\oplus 2}$ as their normal bundles, but which are not  twistor lines.

 Proposition
\ref{prop-nbc2} and Lemmas
\ref{lemma-crtcl-sp} and
\ref{lemma-crtcl-sp2} enable us to determine the normal bundles of
$L_{\theta}^+$ and
$L_{\theta}^-$ in $Z$ for every choices of small resolutions of
$p_{\infty}$. In particular, the normal bundles of $L_{\theta}^+$ and
$L_{\theta}^-$ in $Z$,
and also which component has to be chosen as
 candidates of twistor lines, depend on the choice made.

\subsection{The case of orbit type}\label{ss-orbit} Suppose $\lambda\in I_2\cup I_4$.
In this subsection we calculate the normal bundles of
$L^+_{\alpha}$ and
$L^-_{\alpha}$ in $Z$, where
$L^+_{\alpha}$ and $L^-_{\alpha}$ are  curves which are mapped
biholomorphically onto a real touching conic $C_{\alpha}\in 
\mathcal C_{\lambda}^{\,{\rm{orb}}}$
defined by (\ref{eqn-c}).   Note again that $C_{\alpha}$
and $L^{\pm}_{\alpha}$ depend not only on $\alpha$, but also on
$\lambda\in I_2\cup I_4$. Compared to  generic type and special type,
calculations are much easier since the equations of touching conics of
orbit type are much simpler.

First we make a distinction of
$L^+_{\alpha}$ and 
$L^-_{\alpha}$.  We use  local   coordinates $(x_0,x_1,x_2,z)$
and
$(x_0,x_1,\xi,\eta)$ as in the previous subsection. Recall that 
$\Phi_0^{-1}(H_{\lambda})$ is defined by
$\xi\eta=fx_1^4$, and that  the equation of irreducible components of
$\Phi_0^{-1}(C_{\alpha})$ is  given by $z=\pm
(f-(\alpha+Q)^2)^{\frac{1}{2}}\, x_1^2$ ((\ref{eqn-inv-gen})). Then we
denote by $L^+_{\alpha}$ (resp.
$L^-_{\alpha}$)  the components corresponding to 
$z=(f-(\alpha+Q)^2)^{\frac{1}{2}}\, x_1^2$ (resp.
$z=-(f-(\alpha+Q)^2)^{\frac{1}{2}}\, x_1^2$). $L^+_{\alpha}$ and
$L^-_{\alpha}$ are curves in $Z$.

Recall that in the previous subsection we have introduced an affine
coordinate
$v=\xi/\ell_1\ell_2$  on the exceptional curve $\Gamma_2$. Points on
$\Gamma_2$ are indicated by using this $v$.

\begin{lemma}\label{lemma-intersection3} Let $L^+_{\alpha}$ and
$L^-_{\alpha}$ be as above. Then $ L^+_{\alpha}\cap \Gamma_1,
L^+_{\alpha}\cap\Gamma_3,L^-_{\alpha}\cap\Gamma_1$ and
$L^-_{\alpha}\cap\Gamma_3$ are empty, and $L^+_{\alpha}\cap\Gamma_2$ and
$L^-_{\alpha}\cap\Gamma_2$ are points satisfying respectively
$$ L^+_{\alpha}\cap \Gamma_2=\left\{v=\left(\sqrt{f-(\alpha+Q)^2}+
i(\alpha+Q)\right)\,\frac{x_1^2}{\ell_1\ell_2}\right\}
$$ and
$$ L^-_{\alpha}\cap \Gamma_2=\left\{v=\left(-\sqrt{f-(\alpha+Q)^2}+
i(\alpha+Q)\right)\,\frac{x_1^2}{\ell_1\ell_2}\right\}.
$$ 

\end{lemma}

Here, note again that $\,x_1^2/\ell_1\ell_2$ does not depend on
$x_1$.

\vspace{2mm}
\noindent Proof. Substituting $x_2=\alpha x_1^2$, we have
$$\xi=z+i(x_2+Qx_1^2)=
\left\{\pm \sqrt{f-(\alpha+Q)^2}+i(\alpha+Q)\right\}\,x_1^2$$ and $$\eta=
z-i(x_2+Qx_1^2) =\left\{\pm
\sqrt{f-(\alpha+Q)^2}-i(\alpha+Q)\right\}\,x_1^2$$ over $C_{\alpha}$.
($\pm$ corresponds to $L_{\alpha}^{\pm}$.)
From these and from the explicit resolutions of the
previous subsection, we can easily see  that for any choice of
$\ell_1,\ell_2$ and
$\ell_3$,
$L_{\alpha}^{\pm}\cap \Gamma_1$ and
$L_{\alpha}^{\pm}\cap\Gamma_3$ are empty and that $L_{\alpha}^+\cap
\Gamma_2$ and
$L_{\alpha}^-\cap
\Gamma_2$ are points satisfying 
$$v=\xi/\ell_1\ell_2=\left\{\pm\sqrt{f-(\alpha+Q)^2}+
i(\alpha+Q)\right\}\,\frac{x_1^2}{\ell_1\ell_2},$$ where $\pm$
corresponds to
$L_{\alpha}^+$ and $L_{\alpha}^-$ respectively. Thus we have obtained
all of the claims of the lemma. \proofend

\vspace{2mm} 
Since $L_{\alpha}^{\pm}$ and $\Gamma_2$ are $U(1)$-invariant,
$L_{\alpha}\cap \Gamma_2$ must be  $U(1)$-fixed point.
In particular, any points on $\Gamma_2$ is $U(1)$-fixed.
From these lemmas, we immediately get the following

\begin{lemma}\label{lemma-intersection4} Fix $\lambda\in I_2\cup I_4$.
 Then the set $\{ (L_{\alpha}^+\cup L_{\alpha}^-)\cap\Gamma_2\set
-Q-\sqrt{f}\leq\alpha\leq -Q+\sqrt{f} \}$ is a circle in
$\Gamma_2$ whose center is 
$\Gamma_2\cap
\Gamma_3$ ($=\{v=0\}$) and whose radius is 
$\sqrt{f}\,|x_1^2/\ell_1\ell_2|$. 
\end{lemma}

\vspace{2mm} The following proposition, which corresponds to Propositions
\ref{prop-nottl} (generic type) and
\ref{prop-nbc2} (special type), can be proved by using the same idea as
in  Proposition
\ref{prop-nbc2}. So we omit the proof.

\begin{prop}\label{prop-criet} Set
$h_{2}(\lambda)=\sqrt{f}\,(x_1^2/\ell_1\ell_2)$, which is clearly
differentiable on
 $I_2\cup I_4$. Let $N$ denote the normal bundle of $L^+_{\alpha}$ in
$Z$. Then we have either 
$N\simeq O(1)^{\oplus 2}$ or $N\simeq O\oplus O(2)$, and the latter
holds iff
$\lambda$ is  a critical point of $h_2$.  The same claim holds also for
$L^-_{\alpha}$.
\end{prop}

Needless to say, $h_2$ depends on the choice of $\ell_1$ and $\ell_2$.
Thus as in the case of special type, the normal bundles of $L_{\alpha}^+$
and
$L_{\alpha}^-$ depend on the choice of small resolution of
$p_{\infty}$.  In view of Proposition \ref{prop-criet}, we need to know
the critical point of
$h_2$ for each choice of
$\{\ell_1,\ell_2\}$. There are $4!/(2!2!)=6$  choices of
$\{\ell_1,\ell_2\}$.   If we take $\{\ell_1,\ell_2\}=\{x_0,x_1\}$ for
instance, we have
$h_2(\lambda)^2=(\lambda+1)(a\lambda-b)/\lambda$, and it is elementary 
to determine the critical points of this function. For any other
choices, we always get
$h_2$ in explicit form and it is easy to determine their critical
points.  So here we only present the result:

\begin{lemma}\label{lemma-beha} (i) If $\{\ell_1,\ell_2\}=\{x_0,x_1\}$,
$h_2$ has no critical point on $I_2\cup I_4$. Further,
$h_2(-1)=0,\lim_{\lambda\ua 0}h_2(\lambda)=+\infty,h_2(b/a)=0$ and
$\lim_{\lambda\ua\infty}h_2(\lambda)=+\infty$.  (ii) If
$\{\ell_1,\ell_2\}=\{x_0+x_1,ax_0-bx_1\}$, $h_2$ has no critical point on
$I_2\cup I_4$. Further, $\lim_{\lambda\da
-1}h_2(\lambda)=+\infty,h_2(0)=0,\lim_{\lambda\da
b/a}h_2(\lambda)=+\infty$ and
$\lim_{\lambda\ua\infty}h_2(\lambda)=0$. (iii) If
$\{\ell_1,\ell_2\}=\{x_1,x_0+x_1\}$, $h_2$ has no critical point on
$I_2\cup I_4$. Further, $\lim_{\lambda\da
-1}h_2(\lambda)=+\infty,h_2(0)=0, h_2(b/a)=0$ and
$\lim_{\lambda\ua\infty}h_2(\lambda)=\infty$. (iv) If
$\{\ell_1,\ell_2\}=\{x_0,ax_0-bx_1\}$, $h_2$ has no critical point on
$I_2\cup I_4$. Further,  $h_2(-1)=0,\lim_{\lambda\ua
0}h_2(\lambda)=+\infty,
\lim_{\lambda\da b/a}h_2(\lambda)=+\infty$. (v) If
$\{\ell_1,\ell_2\}=\{x_0,x_0+x_1\}$, or if 
$\{\ell_1,\ell_2\}=\{x_1,ax_0-bx_1\}$, $h_2$ has a unique critical point
on
$I_2$ and
$I_4$ respectively.
\end{lemma}
 

By Corollary \ref{cor-image1}, if $\lambda\in I_2$, images of
twistor lines in
$H_{\lambda}$ must be of orbit type. Therefore by Proposition
\ref{prop-criet},  if a resolution $Z\ra Z_0$ yields a twistor space, 
$h_2$ does not have critical points on $I_2$.  Hence by Lemma
\ref{lemma-beha}, we can conclude that
$\{\ell_1,\ell_2\}\neq\{x_0,x_0+x_1\}$ and 
$\{\ell_1,\ell_2\}\neq\{x_1,ax_0-bx_1\}$. Namely, our investigation
decreases the possibilities of small resolutions.
 We postpone further consequences until the next subsection.

\subsection{Consequences of the results in \S
\ref{ss-generic}--\ref{ss-orbit}}
\label{ss-con} Before stating the results, we again recall our setup.  Let
$B$ be a quartic surface defined by (\ref{eqn-B}) and assume that $Q$ and $f$
 satisfy  the necessary conditions as in Proposition
\ref{prop-necessa}.  Let
$\Phi_0:Z_0\ra\mathbf{CP}^3$ be the double covering branched along $B$.
On $\mathbf{CP}^3$ there is a pencil of $U(1)$-invariant planes
$\{H_{\lambda}\}$, where  $H_{\lambda}$ is defined by $x_0=\lambda x_1$
which is real iff
$\lambda\in\mathbf R\cup\{\infty\}$. For any small resolution $\mu:Z\ra
Z_0$ preserving the real structure, we put $\Phi:=\Phi_0\mu$, and
let
 $\{S_{\lambda}=\Phi^{-1}(H_{\lambda})\}$ be (the real part of) a pencil
of
$U(1)$-invariant divisors on
$Z$, where we put  $S_{\lambda}=\Phi^{-1}(H_{\lambda})$ as before.

We start with the following theorem, which uniquely determines the type
of  real touching conics which can be the images of  twistor lines
contained in
$S_{\lambda}$'s above.

\begin{thm}\label{thm-type}
 Suppose that there is a small resolution  $Z\ra Z_0$ such that $Z$ is a
twistor space.  Let $L$ be a twistor line of $Z$ contained in
$S_{\lambda}$,
$\lambda\in\mathbf R$. Then $\Phi(L)$ is a real touching conic of:  (i) 
special type if
$\lambda\in I_1\cup I_3$, (ii)  orbit type if $\lambda\in I_2$,  (iii)  
generic type if
$\lambda\in I_4$ and if
$\lambda\neq\lambda_0$.
\end{thm}

Note that by Proposition \ref{prop-image}, $\Phi(L)\subset H_{\lambda}$
is a line if $\lambda=\lambda_0$.
 
\vspace{2mm}
\noindent Proof. (i) immediately follows from (ii) of Proposition
\ref{prop-type} and  (v) of Proposition \ref{prop-inv-c}. (ii) is just
Corollary \ref{cor-image1}.
Finally we show   (iii).  By (i) of Proposition \ref{prop-type} it
suffices to show that if
$\lambda\in I_4$, the image cannot be of orbit type. In view of Lemma
\ref{lemma-beha}, we have $h_2(I_2)=(0,\infty)$ and
$h_2(I_4)=(0,\infty)$ for any of the cases (i)--(iv) of the lemma. (We
have already seen that the case (v) can be eliminated.) This implies that  the
circles
 appeared in Lemma
\ref{lemma-intersection3} sweep out
$\Gamma_2\backslash\{\Gamma_2\cap\Gamma_1,\Gamma_2\cap\Gamma_3\}$.
Therefore, $L_{\alpha}^{\pm}\subset S_{\lambda}$ with $\lambda\in
I_2$, and 
$L_{\alpha}^{\pm}\subset S_{\lambda}$ with $\lambda\in I_4$
cannot be the images of twistor lines at the same time. Therefore, if
$\lambda\in I_4$ and if $\lambda\neq \lambda_0$,
the images of twistor lines must be of generic type, as required.
\proofend

\vspace{2mm} The following theorem is the main result of this section.
Recall that 
$p_{\infty}$ is a compound $A_3$-singularity of $Z_0$, and there are
$4!=24$ choices of  small resolutions of $p_{\infty}$
(see \S\ref{ss-special}). Recall also that
once a resolution of
$p_{\infty}$ is given, it naturally induces that of $\ol{p}_{\infty}$ by
reality.

\begin{thm}\label{thm-sr} Among 24 ways of possible small resolutions of
$p_{\infty}$, 22 resolutions do not yield a twistor space. The
remaining two  resolutions are given by the following two choices of
linear forms:
$$\ell_1=x_1,\hspace{2mm}\ell_2=x_0+x_1,\hspace{2mm}\ell_3=x_0,$$ and
$$\ell_1=ax_0-bx_1,\hspace{2mm}\ell_2=x_0,\hspace{2mm}\ell_3=x_0+x_1.$$
\end{thm}

Here we do not yet claim that the threefolds obtained by these two
resolutions are actually twistor spaces.

\vspace{2mm}
\noindent Proof. By Theorem \ref{thm-type} (i), if $\lambda\in I_1$, the
images of twistor lines in
$S_{\lambda}$ are real touching conics  of special type. As in Section
\ref{ss-special}, there are two families $\mathcal L_{\lambda}^+$ and
$\mathcal L_{\lambda}^-$ of real rational curves which are candidates of
 twistor lines in $S_{\lambda}$. As we have already remarked in
Section \ref{s-inv},
$L^+_{\theta}\in\mathcal L_{\lambda}^+$ and 
$L^-_{\theta}\in\mathcal L_{\lambda}^-$ 
cannot be twistor lines simultaneously. Suppose first that (any of the)
members of $\mathcal L_{\lambda}^+$  are twistor lines. Then by
Proposition \ref{prop-nbc2}, the function $h_1$ does not have critical
points on $I_1$. By Lemma
\ref{lemma-crtcl-sp}, this implies that  we have either
\begin{equation}\label{eqn-candi1}
\ell_1=x_1 \hspace{3mm}{\rm{or}}\hspace{3mm}
\ell_1=x_0+x_1.\end{equation} On the other hand, by Corollary
\ref{cor-oneof} (i), under our assumption, members of $\mathcal
L_{\lambda}^-$ are twistor lines for $\lambda\in I_3$. Therefore again
by Proposition
\ref{prop-nbc2},
$h_3$ does not have critical points on
$I_3$. Then by Lemma \ref{lemma-crtcl-sp2}, the cases (i) and (iii)
of the lemma are
eliminated and  we have either 
\begin{equation}\label{eqn-candi2}
\{\ell_1,\ell_2,\ell_3\}=\{x_1,x_0+x_1,ax_0-bx_1\}
\hspace{3mm}{\rm{or}}\hspace{3mm}
\{\ell_1,\ell_2,\ell_3\}=\{x_0,x_1,x_0+x_1\}.
\end{equation}(Note here that we do not specify the order.) 

Next we consider twistor lines in $S_{\lambda}$ for $\lambda\in I_2$. By 
Theorem \ref{thm-type} (ii) the images are real touching conics of orbit
type. Then  Proposition \ref{prop-criet} implies that $h_2$ has no
critical point on
$I_2$. Hence by Lemma \ref{lemma-beha},  we have either 
\begin{equation}\label{eqn-candi3}
\{\ell_1,\ell_2\}=\{x_0,x_1\}
\hspace{2mm}{\rm{or}}\hspace{2mm}
\{x_0+x_1,ax_0-bx_1\}
\hspace{2mm}{\rm{or}}\hspace{2mm}
\{x_1,x_0+x_1\}
\hspace{2mm}{\rm{or}}\hspace{2mm}
\{x_0,ax_0-bx_1\}.
\end{equation} 

 Now we note other
restrictions: namely,  when
$\lambda$ increases to pass from $I_1$ to
$I_2$,  twistor lines in $S_{\lambda}$ varies continuously, so that  we
have 
\begin{equation}\label{eqn-limit-1}
\lim_{\lambda\ua-1}h_1(\lambda)=\left(\lim_{\lambda\da-1}h_2(\lambda)\right)^{-1}.
\end{equation} (Here the inverse of the right hand side is a consequence
of the fact that 
$\Gamma_1\cap \Gamma_2=\{u=0\}=\{v=\infty\}$.) Similarly, moving
$\lambda$ from
$I_2$ to $I_3$, we have 
\begin{equation}\label{eqn-limit-2}
\lim_{\lambda\ua 0}h_2(\lambda)=\left(\lim_{\lambda\da
0}h_3(\lambda)\right)^{-1}.
\end{equation} Take $\ell_1=x_1$ for the first example. Then by Lemma
\ref{lemma-crtcl-sp} (i) we have
$h_1(-1)=0$. Hence it follows from (\ref{eqn-limit-1}) that  
$\lim_{\lambda\da-1}h_2(\lambda)=\infty$. Then the cases (i) and (iv)
of Lemma \ref{lemma-beha} fail and we have 
$\{\ell_1,\ell_2\}=\{x_0+x_1,ax_0-bx_1\}$ ((ii)) or
$\{\ell_1,\ell_2\}=\{x_1,x_0+x_1\}$ ((iii)). The former clearly fails and
we get
$\ell_2=x_0+x_1$. This appears  in (\ref{eqn-candi3}). Then we have from
Lemma
\ref{lemma-beha} (iii) that $h_2(0)=0$. Hence by (\ref{eqn-limit-2}), we
have
$\lim_{\lambda\da 0}h_3(\lambda)=\infty$. It then follows from Lemma
\ref{lemma-crtcl-sp2} that $\ell_3=x_0$. Thus we get
$\ell_1=x_1,\ell_2=x_0+x_1,\ell_3=x_0$.

Next take $\ell_1=x_0+x_1$. Then we have $\lim_{\lambda\uparrow
-1}h_1(\lambda)=+\infty$ (Lemma \ref{lemma-crtcl-sp} (iii)), so that
$\lim_{\lambda\da-1}h_2(\lambda)=0$. Then looking (i)--(iv) of Lemma
\ref{lemma-beha}, this possibility fails. Namely, 
we have $l_1\neq x_0+x_1$. Thus we can conclude that if
$L^+_{\theta}\in\mathcal L^+_{\lambda}$ is a twistor line over $I_1$,  it
follows that
$\ell_1=x_1,\ell_2=x_0+x_1$ and $\ell_3=x_0$.  This is the former
candidate of the theorem. 

Next suppose  that  $L^-_{\theta}\in\mathcal L^-_{\lambda}$ is a twistor
line over $I_1$ and repeat similar argument above. By Proposition
\ref{prop-nbc2}, $h_3$ has no critical point on $I_1$. It then follows
from Lemma
\ref{lemma-crtcl-sp2} that  either
$\{\ell_1,\ell_2,\ell_3\}=\{x_0,x_0+x_1,ax_0-bx_1\}$ ((i)) or
$\{\ell_1,\ell_2,\ell_3\}=\{x_0,x_1,ax_0-bx_1\}$ ((iii)) holds.
On the other hand,  (\ref{eqn-candi3}) is valid also in this case.
Further we have as before
\begin{equation}\label{eqn-fsfa}
\lim_{\lambda\ua-1}h_3(\lambda)=\left(\lim_{\lambda\da-1}h_2(\lambda)\right)^{-1}
\hspace{2mm}{\rm{and}}\hspace{3mm}
\lim_{\lambda\ua 0}h_2(\lambda)=\left(\lim_{\lambda\da
0}h_1(\lambda)\right)^{-1}.
\end{equation}

If $\{\ell_1,\ell_2,\ell_3\}=\{x_0,x_0+x_1,ax_0-bx_1\}$, then
$\lim_{\lambda\ua -1}h_3(\lambda)=+\infty$ (Lemma \ref{lemma-crtcl-sp2}
(i)), so that we have 
$\lim_{\lambda\da -1}h_2(\lambda)=0$ by (\ref{eqn-fsfa}). Hence by Lemma
\ref{lemma-beha} we have 
$\{\ell_1,\ell_2\}=\{x_0,ax_0-bx_1\}$, which implies 
$\lim_{\lambda\ua 0}h_2(\lambda)=\infty$ ((iv) of Lemma
\ref{lemma-beha}) and
$\ell_3=x_0+x_1$. Hence $\lim_{\lambda\da 0}h_1(\lambda)=0$
by (\ref{eqn-fsfa}).  It follows
from Lemma \ref{lemma-crtcl-sp} that $\ell_1=ax_0-bx_1$, which means
$\ell_2=x_0$.   (iv) of Lemma \ref{lemma-crtcl-sp} says that $h_1$ has
no critical point on
$I_3$, which is consistent with the fact that 
$L^+_{\theta}$ is a twistor line over $I_3$.

If $\{\ell_1,\ell_2,\ell_3\}=\{x_0,x_1,ax_0-bx_1\}$, 
$h_3(-1)=0$ (Lemma \ref{lemma-crtcl-sp2} (iii)), so that we have 
$\lim_{\lambda\da -1}h_2(\lambda)=\infty$ by (\ref{eqn-fsfa}). Therefore
we get the two possibilities (ii) and (iii) of  Lemma
\ref{lemma-beha}, but both contain $x_0+x_1$ which is not compatible
with our choice of $\{\ell_1,\ell_2,\ell_3\}$. Thus we have
$\{\ell_1,\ell_2,\ell_3\}\neq \{x_0,x_1,ax_0-bx_1\}$.
 This implies that  if $L^-_{\theta}\in\mathcal L^-_{\lambda}$ is a
twistor line for $\lambda\in I_1$,  then 
$\ell_1=ax_0-bx_1,\ell_2=x_0$, and $\ell_3=x_0+x_1$.  This is the latter
candidate of the theorem, and we have completed the proof.
\proofend

\vspace{2mm} At first sight it may not be evident why there are two
choices of  small resolutions which can yield twistor spaces. But our  
proof shows the difference of them.  To explain this,
for each
$\lambda\in I_1\cup I_2\cup I_3$, let
$L_{\lambda}\subset Z$ be any  of the members of $\mathcal
L^{\pm}_{\lambda}$ which are chosen as  candidates  of twistor lines in
the proof of Theorem
\ref{thm-sr}. Namely, if 
$\ell_1=x_1,\ell_2=x_0+x_1,\ell_3=x_0$ (the former case),  members of
$\mathcal L_{\lambda}^+$ (resp. $\mathcal L_{\lambda}^-$) must  be
chosen for
$\lambda\in I_1$ (resp.
\!$\lambda\in I_3$).  If $\ell_1=ax_0-bx_1,\ell_2=x_0,\ell_3=x_0+x_1$
(the latter case),  members of $\mathcal L_{\lambda}^-$ (resp. $\mathcal
L_{\lambda}^+$) have to be chosen for $\lambda\in I_1$ (resp.
$\lambda\in I_3$). (For $\lambda\in I_2$   any
members of $\mathcal L_{\lambda}^{+}$ and $\mathcal L_{\lambda}^{-}$
must be chosen simultaneously as in Proposition
\ref{prop-inv-c}.) Our proof shows that  when $\lambda\in \mathbf R$ 
increases from $-\infty$ to
$b/a$,  the intersection $\Gamma\cap L_{\lambda}$ moves from $\Gamma_1$
to
$\Gamma_3$ for the former choice, whereas it moves from $\Gamma_3$ to
$\Gamma_1$ for the latter choice. Namely,  exchanging the choice
reverses the orientation of $\Gamma\cap L_{\lambda}$ as $\lambda$
increases.

By  similar consideration  we can determine which irreducible component
must be chosen over
$I_4^-\cup I_4^+$.  First take the small resolution of $p_{\infty}$
determined by the first choice of $\ell_1,\ell_2$ and $\ell_3$ in
Theorem \ref{thm-sr}. Then as mentioned above, over $I_3$, members of
$\mathcal L^-_{\lambda}$ must be chosen. Further, we
have
$h_3(b/a)=0$ by Lemma \ref{lemma-crtcl-sp2} (iv). This implies that
as $\lambda$ goes to $b/a$,
the circle of intersection
$\cup\{\Gamma_3\cap L^-_{\theta}\set
 L^-_{\theta}\in\mathcal L^-_{\lambda}\}$ shrinks
to be the $U(1)$-fixed point of
$\Gamma_3$ which is different from $\Gamma_2\cap \Gamma_3$.
 On the other hand, on the irreducible component of
$\Phi^{-1}(\ell_{\infty})$ which intersects
 $\Gamma_3$, and which is different from $\Gamma_2$, one can use
$x_2=y_2/y_3$ as an affine coordinate whose center is the intersection
point with $\Gamma_3$.  Therefore,  the circle of intersection
 appeared in   Lemma \ref{lemma-int-a} also must shrink to be the origin as
$\lambda$ decreases to be $b/a$.  This uniquely determines which one of $\mathcal
L_{\lambda}^+$ and $\mathcal L_{\lambda}^-$  has to  be chosen for
$\lambda\in I_4^-$. Then  for
$\lambda\in I_4^+$ another irreducible component must be chosen. 
The case of the latter choice
of
$\ell_1,\ell_2$ and
$\ell_3$ is now obvious. Namely,  we have to choose the irreducible component
over
$I_4^-$ and
$ I_4^+$ which are different from the former case.

\section{Twistor lines whose images are lines}\label{s-line}

In this section we study real lines whose images are lines
in $\mathbf{CP}^3$.
Recall  we have shown in Proposition \ref{prop-image} that if
$L$ is a real line  intersecting $\Gamma_0$, then $\Phi(L)$ is a line
 going through $P_0$. The following proposition is its
 converse.

\begin{prop}\label{prop-imline} Let $l\subset\mathbf{CP}^3$ be any real line
going through
$P_0$. Then $\Phi^{-1}(l)$ has just two irreducible components, both of which are
real, smooth and rational. One of the components is the exceptional
curve $\Gamma_0$ and  another component is  mapped (2 to 1) onto $l$.
Further, the normal bundle of the latter component in $Z$ is isomorphic
to $O(1)^{\oplus 2}$. 
\end{prop}

\noindent Proof. First we note that if $l$ is a  real line,  $B\cap
l$ consists of just three points, one of which is
$P_0$. This follows from the facts that, $B$ is a quartic,  $B\cap l$ is
real,
$P_0$ is the unique real point of $B$ (Proposition
\ref{prop-realpoint}), and $P_0$ is a double point. Therefore
$\Phi_0^{-1}(l)\ra l$ is  two-to-one covering branched at
three points. Let $P$ and
$\ol{P}$ be the two branch points other than $P_0$. 
Because $l$ intersects $B$ transversally at these two points,
$\Phi_0^{-1}(P)$ and $\Phi_0^{-1}(\ol{P})$ are smooth points of
$\Phi_0^{-1}(l)$. Further, since $P_0$ is an ordinary double point of $B$,
$\Phi_0^{-1}(P_0)$ is a node of
$\Phi_0^{-1}(l)$. From these it  follows that  
$\Phi^{-1}(l)$ has just two irreducible components, one of which is
$\Gamma_0$. Let
$L$ be the irreducible component other than $\Gamma_0$. Then $L$ is 
smooth and
$L\ra l$ is two-to-one covering whose branch points are $P$ and
$\ol{P}$.  
 Therefore by Hurwitz, $L$ is a  rational curve.
$L$ is  real since $\Phi_0^{-1}(l)$ is real.

It remains to show that $N_{L/Z}\simeq  O(1)^{\oplus 2}$. The idea is
similar to   Propositions \ref{prop-nbc},
\ref{prop-nbc2}, and
\ref{prop-criet}. We first show that $N_{L/Z}\simeq O(1)^{\oplus 2}$
 or $N_{L/Z}\simeq O\oplus O(2)$. By Bertini, $H\cap B$ is smooth
outside $P_0$ for general plane $H$ containing $ l$. Further, since
$H\cap B$ is a quartic,
$S:=\Phi^{-1}(H)$ is a smooth rational surface with $c^2_1=2$. Moreover,
$\Phi^{-1}(l)$ is an anticanonical curve of
$S$ so that we have $(\Gamma_0+L)^2=2$ on $S$. Furthermore, it is readily
seen that
$\Gamma_0^2=-2$ and $\Gamma_0\cdot L=2$ on $S$. Therefore we have $L^2=0$ on $S$.
Then the argument in the proof of Proposition \ref{prop-nbc} implies 
$N_{L/Z}\simeq O(1)^{\oplus 2}$ or $N_{L/Z}\simeq O(2)\oplus O$.
 To show that the latter does
not hold, we first see that $\Gamma_0/\langle\sigma\rangle$ is canonically
identified with the projective space of real lines going through $P_0$. Concretely,
for each real line $l\ni P_0$, we associate the intersection $\Gamma_0\cap
(\Phi^{-1}(l)-\Gamma_0)$ which is a conjugate pair of points.
We show by explicit calculation that this correspondence,
which we will denote by $\psi$, is actually
an isomorphism. The problem being local, we use a local coordinate
$(w_1,w_2,w_3)$  (around $P_0$) defined in (\ref{coord}).  Then in a
neighborhood of
$P_0$, $Z_0$ is given by the equation 
\begin{equation}\label{eqn-Z_0}z^2+w_1^2+w_2w_3=0.\end{equation}  Small
resolutions of the double point $p_0=\Phi_0^{-1}(P_0)\in Z_0$ are
explicitly obtained by blowing-up along $\{z+iw_1=w_3=0\}$ or
$\{z+iw_1=w_2=0\}$. In the former case, we can use
$(z+iw_1:w_3)=(-w_2:z-iw_1)$ as a homogeneous coordinate on $\Gamma_0$,
whereas in the latter case we can use
$(z+iw_1:w_2)=(-w_3:z-iw_1)$ instead. We see only in the former case,
since the calculation is identical. Let
$(w_1:w_2:w_3)$ be a real line through $P_0$.  Namely, we assume
$w_1\in\mathbf R$,  $\ol{w}_2=w_3$, and $w_1^2+|w_2|^2\neq 0$.  Then by
(\ref{eqn-Z_0}), we have $z=\pm i(w_1^2+|w_2|^2)^{1/2}$. Hence we get
$(z+iw_1:w_3)=(i(w_1\pm (w_1^2+|w_2|^2)^{1/2}):\ol{w}_2)$. Namely,
$\psi$ is explicitly given by 
\begin{equation}\label{eqn-correspon} 
\psi:(w_1:w_2:\ol{w}_2)\longmapsto 
\left(i\left(w_1\pm \sqrt{w_1^2+|w_2|^2}\right):\ol{w}_2\right).
\end{equation} 
(Note that the image of (\ref{eqn-correspon}) is considered as a point of 
$\Gamma_0/\langle\sigma\rangle$.)
First suppose 
$w_1\neq 0$. It is readily seen that we can suppose $w_1=1$. Then  in
(\ref{eqn-correspon})  the image becomes 
$(i(1\pm(1+|w_2|^2)^{1/2}):\ol{w}_2)$. Taking the sign `$+$', 
(\ref{eqn-correspon}) can be rewritten as
\begin{equation}\label{eqn-correspon2}
\psi:\mathbf C\ni w_2\mapsto 
\frac{-i\ol{w}_2}{1+\sqrt{1+|w_2|^2}},
\end{equation} where we use (the second entry)$/$(the first entry) as
an affine coordinate on $\Gamma_0$. The image of (\ref{eqn-correspon2}) is
clearly contained in the unit disk
$\{u\in\mathbf C\set |u|<1\}$. We show that (\ref{eqn-correspon2})
give a diffeomorphism between
$\mathbf C=\mathbf R^2$ and the unit disk. Putting $w_2=re^{i\theta}$,
(\ref{eqn-correspon2})  is rewritten as
\begin{equation}\label{eqn-correspon3}
\psi: re^{i\theta}\longmapsto
\frac{-ire^{-i\theta}}{1+\sqrt{1+r^2}}.
\end{equation} It is elementary to show that
$k(r):=r/(1+\sqrt{1+r^2})$ is differentiable on
$\{r>0\}$ and  its derivative is always positive, and that
$\lim_{r\ua\infty}k(r)=1$ and $\lim_{r\da 0}k(r)=0$ hold.
Hence $k$ gives  a bijection between $\{r\geq 0\}$ and
$\{0\leq s<1\}$.
It follows that (\ref{eqn-correspon2}) gives a bijection between
$\mathbf C$ and the unit disk.
 Moreover, the positivity of $k'$ implies 
that (\ref{eqn-correspon3}) is a diffeomorphism on $\mathbf C^*$.
For $w_2=0$, it can be easily checked that 
$(\partial w_2/\partial \ol{w}_2)(0)\neq 0$. Therefore (\ref{eqn-correspon2})
is a diffeomorphism on $\mathbf C$.

Next consider the case $w_1=0$. Then  we have $w_2\neq 0$, and the
image becomes
$(\pm i|w_2|:\ol{w}_2)=(1:\pm i\ol{w}_2/|w_2|)$. From this, it easily
follows that (\ref{eqn-correspon}) gives a diffeomorphism between the two
subsets
$\{(0:1:w)\set w\in U(1)\}$ and 
$\{(1:u)\in \Gamma_0\set u\in U(1)\}/\langle\sigma\rangle$. 
Moreover on $\mathbf{RP}^2\backslash\mathbf R^2$ we can use
$(1/r,\theta)$ as a local coordinate on $\mathbf{RP}^2$.
Then we can readily show that $(d/ds)(k(1/s))|_{s=0}\neq 0$.
This implies that $\psi$ is diffeomorphic also in a neighborhood of
$\mathbf{RP}^2\backslash\mathbf R^2$, the infinite circle.
Note that the bijectivity of $\psi$ 
 implies that if
$l\neq l'$, then the corresponding rational curves in $Z$ are
disjoint. 

Next take any  real plane $H$ containing $l$.
On $H$ there is a one-dimensional family of lines through
$P_0$. Taking the inverse image, we obtain a one-dimensional 
holomorphic family $\mathcal L_{H}$
of  rational curves in $Z$, containing $L$ as a real member.
Any real member of $\mathcal L_{H}$ defines a conjugate pair of points
as the intersection with $\Gamma_0$. 
Consequently, $\mathcal L_H$ determines a real circle $\mathcal C_H$
 in $\Gamma_0$.
(Namely $\mathcal C_H=\{L'\cap \Gamma_0\set L'\in\mathcal L_H^{\sigma}\}$.)
If $s$ denotes the section of $N=N_{L/Z}$ associated to $\mathcal L_H$,
then Re$s(z)$ is non-zero by the diffeomorphicity of $\psi$, and
is represented by a tangent vector of
$\mathcal C_{H}$ at $z$, where we put $\{z,\ol{z}\}=\Gamma_0\cap L$.

Let $\{v_1,v_2\}$ be any oriented orthogonal basis of $T_z\Gamma_0$,
where we take the complex orientation and orthogonality.
Then since we have the isomorphism $\psi$,
there is a unique real plane $H_i$ ($i=1,2$) containing $l$ such that
$v_i$ is tangent to $\mathcal C_{H_i}$.
Let $s$ (resp.\,$t$)   be the global section of $N=N_{L/Z}$
associated to $\mathcal L_{H_1}$ (resp.\,$\mathcal L_{H_2}$).
 We now claim that
$as+bt$ does not vanish at $z$ and $\ol{z}$ simultaneously, 
unless $(a,b)= (0,0)$.
Putting $a=a_1+ia_2$ and $b=b_1+ib_2$, we readily have
\begin{equation}\label{eqn-reim4}
{\rm{Re}}(as+bt)=(a_1{\rm{Re}}s-b_2{\rm{Im}}t)+(b_1{\rm{Re}}t-a_2{\rm{Im}}s).
\end{equation}
Since (Re$s)(z)$ and (Re$s)(\ol{z})$ 
(resp.\,(Re$t)(z)$ and (Re$t)(\ol{z})$) are represented by  tangent vectors of
$\mathcal C_{H_1}$ (resp. $\mathcal C_{H_2}$), our choice of $H_1$ and $H_2$
implies that $({\rm{Re}}s)(z)$ is parallel to $({\rm{Im}}t)(z)$ and
$({\rm{Re}}t)(z)$ is parallel to $({\rm{Im}}s)(z)$.
The same is true at $\ol{z}$.
Hence by (\ref{eqn-reim4})
 if ${\rm{Re}}(as+bt)(z)=0$, then $a_1{\rm{Re}}s(z)=b_2{\rm{Im}}t(z)$
and $b_1{\rm{Re}}t(z)=a_2{\rm{Im}}s(z)$,
and ${\rm{Re}}(as+bt)(\ol{z})=0$ implies similar equalities.
Since ${\rm{Re}}s, {\rm{Re}}t,{\rm{Im}}s$ and ${\rm{Im}}t$ do not
be zero at both of $z$ and $\ol{z}$  as is already mentioned,
$a_1=0$ iff $b_2=0$ and $a_2=0$ iff $b_1=0$.
Therefore either $a_1b_2\neq 0$ or $b_1a_2\neq 0$ holds.
Suppose $a_1b_2\neq 0$. Then we show that 
$a_1{\rm{Re}}s(z)=b_2{\rm{Im}}t(z)$ and 
$a_1{\rm{Re}}s(\ol{z})=b_2{\rm{Im}}t(\ol{z})$ cannot hold simultaneously:
suppose that $a_1b_2>0$. Then 
${\rm{Re}}s(z)$ and ${\rm{Im}}t(z)$ have the same direction and
it follows that $\{{\rm{Re}}t(z),{\rm{Re}}s(z)\,
(=(b_2/a_1){{\rm{Im}}t(z)})\} $ is an
oriented basis of
$T_z\Gamma_0$. On the other hand, we have  
 ${\rm{Re}}s(\ol{z})=\sigma_*({\rm{Re}}s(z))$ and
 ${\rm{Re}}t(\ol{z})=\sigma_*({\rm{Re}}t(z))$. 
Therefore $\{{\rm{Re}}t(\ol{z}), {\rm{Re}}s(\ol{z})\}$ is an anti-oriented
basis of $T_{\ol{z}}\Gamma_0$ because $\sigma$ is orientation reversing.
On the other hand, $a_1{\rm{Re}}s(\ol{z})=b_2{\rm{Im}}t(\ol{z})$ and
$a_1b_2\neq 0$ imply that 
$\{{\rm{Re}}t(\ol{z}),{\rm{Re}}s(\ol{z})\} $ is an oriented basis of 
 $T_{\ol{z}}\Gamma_0$. This is a contradiction.
The case $b_1a_2>0$ is similar.
Therefore ${\rm{Re}}(as+bt)$ cannot be zero at 
$z$ and $\ol{z}$ simultaneously  provided $a_1b_2\neq 0$.
If $b_1a_2\neq 0$, then $b_1{\rm{Re}}t(z)=a_2{\rm{Im}}s(z)$ and 
$b_1{\rm{Re}}t(\ol{z})=a_2{\rm{Im}}s(\ol{z})$ do not hold at the same time.
 Thus we have shown that ${\rm{Re}}(as+bt)$ cannot be zero at 
$z$ and $\ol{z}$ simultaneously for any $(a,b)\neq (0,0)$. 
Hence so does $as+bt$. Therefore we get
$N\simeq O(1)^{\oplus 2}$ by Lemma \ref{lemma-nb}.
\proofend

\vspace{2mm} Thus in our complex manifold $Z$ there actually exists a
connected family  of real rational curves parametrized by
$\Gamma_0/\langle\sigma\rangle\simeq\mathbf{RP}^2$, whose normal bundle is
isomorphic to $O(1)^{\oplus 2}$. 
By Proposition \ref{prop-image},  and the canonical isomorphism 
$\psi$ obtained in the previous proof,
all of these real lines must be twistor lines (if $Z$ is a twistor space).
 Obviously this family
 is  $U(1)$-invariant, although
general members are not 
$U(1)$-invariant:

\begin{prop} Among this family of real lines in $Z$, just one member is
$U(1)$-invariant. Further, the member is fixed  by $U(1)$
pointwisely.
\end{prop}

\noindent Proof. 
Recall that in a neighborhood of $p_0$,
$Z_0$ is defined by the equation
$z^2+w_1^2+w_2w_3=0$ ((\ref{eqn-Z_0})). It is immediate to see that the
$U(1)$-action looks like $(w_1,w_2,w_3)\mapsto
(w_1,tw_2,t^{-1}w_3)$ for
$t\in U(1)$. Thus using  homogeneous coordinates used in the last
proof, the $U(1)$-action on $\Gamma_0$ is given by 
$(u:v)\mapsto (u:tv)$ or $(u:v)\mapsto (u:t^{-1}v)$, depending on the
choice of a small resolution of $p_0$. Therefore only the real line
corresponding to 
$[(1:0)](=[(0:1)])\in \Gamma_0/\langle\sigma\rangle$ is $U(1)$-fixed. 
In view of (\ref{eqn-correspon}) and (\ref{coord}), the equation of this
line is explicitly given by
$y_2=y_3=0$, which is   pointwisely  $U(1)$-fixed by Proposition \ref{prop-def-B}.
 Since
$\Phi:Z\ra\mathbf{CP}^3$ is 
$U(1)$-equivariant, it follows that the corresponding rational
curve in
$Z$ is also pointwisely fixed. \proofend

\vspace{2mm}
What we have done in this paper  can be summarized as follows:
\begin{prop}
Let $Z_0$ be as in Proposition \ref{prop-def-B}, where $Q$ and $f$
satisfy the conditions in Proposition \ref{prop-necessa} which is
necessary for $Z_0$ to be birational to a twistor space. 
Let $\mu:Z\ra Z_0$ be one of the small resolutions determined in Theorem
\ref{thm-sr},
where the real ordinary double point is resolved in arbitrary way (as a small resolution).
Set $\Phi=\mu\Phi_0$ and consider
the curve $\Phi^{-1}(l_{\infty})$  which is a real cycle of $U(1)$-invariant smooth rational curves  consisting of eight
irreducible components
Then for any smooth points of $\Phi^{-1}(l_{\infty})$,
we can explicitly specify a real smooth rational curve in $Z$ going
through the points whose normal bundle is $ O(1)^{\oplus 2}$. 
Furthermore, if $Z$ is a twistor space, all  the curves we specified
must be  twistor lines.
\end{prop}

\small
\vspace{13mm}
\hspace{4.5cm}
$\begin{array}{l}
\mbox{Department of Mathematics}\\
\mbox{Graduate School of Science and Engineering}\\
\mbox{Tokyo Institute of Technology}\\
\mbox{2-12-1, O-okayama, Meguro, 152-8551, JAPAN}\\
\mbox{{\tt {honda@math.titech.ac.jp}}}
\end{array}$

\end{document}